\input amsppt.sty
\magnification=\magstephalf \hsize = 6.5 truein \vsize = 9 truein
\vskip 3.5 in

\NoBlackBoxes
\TagsAsMath

\def\label#1{\par%
        \hangafter 1%
        \hangindent .75 in%
        \noindent%
        \hbox to .75 in{#1\hfill}%
        \ignorespaces%
        }

\newskip\sectionskipamount
\sectionskipamount = 24pt plus 8pt minus 8pt
\def\sectionskip{\vskip\sectionskipamount}
\define\sectionbreak{%
        \par  \ifdim\lastskip<\sectionskipamount
        \removelastskip  \penalty-2000  \sectionskip  \fi}
\define\section#1{%
        \sectionbreak   
        \subheading{#1}%
        \bigskip
        }

\redefine\qed{{\unskip\nobreak\hfil\penalty50\hskip2em\vadjust{}\nobreak\hfil
    $\square$\parfillskip=0pt\finalhyphendemerits=0\par}}

        
        \let    \< = \langle
        \let    \> = \rangle

\define\op#1{\operatorname{\fam=0\tenrm{#1}}} 
        \define         \a              {\alpha}
        \redefine       \b              {\beta}
        \redefine       \d              {\delta}
        \redefine       \D              {\Delta}
        \define         \e              {\varepsilon}
        \define         \E              {\op {E}}
        \define         \g              {\gamma}
        \define         \G              {\Gamma}
        \redefine       \l              {\lambda}
        \redefine       \L              {\Lambda}
        \define         \n              {\nabla}
        \redefine       \var            {\varphi}
        \define         \s              {\sigma}
        \redefine       \Sig            {\Sigma}
        \redefine       \t              {\tau}
        \define         \th             {\theta}
        \redefine       \O              {\Omega}
        \redefine       \o              {\omega}
        \define         \z              {\zeta}
        \define         \k              {\kappa}
        \redefine       \i              {\infty}
        \define         \p              {\partial}
        \define         \vsfg           {\midspace{0.1 truein}}

\document

\topmatter
\title On certain $C$-test words for free groups
\endtitle
\author  Donghi Lee \endauthor

\address {Department of Mathematics, University of Illinois at
Urbana--Champaign, 1409 West Green Street, Urbana, IL 61801, USA}
\endaddress

\email {d-lee9\@math.uiuc.edu}
\endemail


\abstract{Let $F_m$ be a free group of a finite rank $m \ge 2$ and
$X_i, \, Y_j$ be elements in $F_m$. A non-empty word $w(x_1,
\dots, x_n)$ is called a $C$-test word in $n$ letters for $F_m$ if,
whenever $w(X_1, \dots, X_n)=w(Y_1, \dots, Y_n) \neq 1$, the two
$n$-tuples $(X_1, \dots, X_n)$ and $(Y_1, \dots, Y_n)$ are
conjugate in $F_m$. In this paper we construct, for each $n \ge
2$, a $C$-test word $v_n(x_1, \dots, x_n)$ with the additional
property that $v_n(X_1, \dots, X_n)=1$ if and only if the
subgroup of $F_m$ generated by $X_1, \dots, X_n$ is cyclic.
Making use of such words $v_m(x_1, \dots, x_m)$ and $v_{m+1}(x_1,
\dots, x_{m+1})$, we provide a positive solution to the following
problem raised by Shpilrain: There exist two elements $u_1, u_2
\in F_m$ such that every endomorphism $\psi$ of $F_m$ with
non-cyclic image is completely determined by $\psi (u_1)$, $\psi
(u_2)$.}
\endabstract
\endtopmatter


\baselineskip=24pt

\heading 1. Introduction
\endheading

Let $F_m=\<x_1, \dots, x_m\>$ be the free group of a finite rank $m
\ge 2$ on the set $\{x_1, \dots, x_m\}$. The purpose of this
paper is to present a positive solution to the following problem
raised by Shpilrain [1]:

\noindent {\bf Problem.} Are there 2 elements $u_1$, $u_2$ in
$F_m$ such that any endomorphism $\psi$ of $F_m$ with non-cyclic
image is uniquely determined by $\psi(u_1)$, $\psi(u_2)$ ? (In
other words, are there 2 elements $u_1$, $u_2$ in $F_m$ such that
whenever $\phi(u_i)=\psi(u_i)$, $i=1,\, 2$, for endomorphisms $\phi$, $\psi$ of $F_m$ with non-cyclic images, it follows that
$\phi=\psi$ ?)

In [2], Ivanov solved in the affirmative this
problem in the case where $\psi$ is a monomorphism of $F_m$
by constructing a so-called $C$-test word $w_n(x_1, \dots, x_n)$ for
each $n \ge 2$. A $C$-test word is defined due to Ivanov [2] as
follows:

\noindent {\bf Definition.} A non-empty word $v(x_1, \dots, x_n)$
is a {\it $C$-test word} in $n$ letters for $F_m$ if for any two
$n$-tuples $(X_1, \dots, X_n)$, $(Y_1, \dots, Y_n)$ of elements
of $F_m$ the equality $v(X_1, \dots, X_n) = v(Y_1, \dots, Y_n)
\neq 1$ implies the existence of an element $S \in F_m$ such
that $Y_i = S X_i S^{-1}$ for all $i = 1,\, 2,\, \dots,\, n$.

According to the result of [2, Corollary 1], if $v$ is a $C$-test
word in $m$ letters for $F_m$, $\phi$ is an endomorphism, $\psi$
is a monomorphism of $F_m$, and $\phi(v)=\psi(v)$, then we have
$\phi=\tau_S \circ \psi$, where $S \in F_m$ is such that $\<S,
\psi(v)\>$ is cyclic, and $\tau_S$ is the inner automorphism of
$F_m$ defined by means of $S$. Notice that this assertion is no
longer true if $\psi$ is extended to an endomorphism of $F_m$
with non-cyclic image, since the fact $\psi(v) \neq 1$ is no
longer guaranteed. The efforts to extend the result of [2,
Corollary 1] to the case where $\psi$ is an endomorphism of $F_m$
with non-cyclic image have led us to proving the following:

\proclaim {Theorem} For every $n \ge 2$ there exists a $C$-test
word $v_n(x_1, \dots, x_n)$ in $n$ letters for $F_m$ with the
additional property that $v_n(X_1, \dots , X_n)=1$ if and only if
the subgroup $\< X_1, \dots, X_n \>$ of $F_m$ generated by $X_1,
\dots, X_n$ is cyclic.
\endproclaim

We construct such a $C$-test word $v_n(x_1, \dots, x_n)$ by combining
Ivanov's $C$-test word $w_2(x_1, x_2)$ and an auxiliary word
$u(x_1, x_2)$ defined below. Here, let us recall Ivanov's
$C$-test word $w_2(x_1, x_2)$:
$$
\multline w_2(x_1,x_2) = \,[x_1^8,x_2^8]^{100} x_1
[x_1^8,x_2^8]^{200}
x_1 [x_1^8,x_2^8]^{300} x_1^{-1} [x_1^8,x_2^8]^{400} x_1^{-1} \\
[x_1^8,x_2^8]^{500} x_2 [x_1^8,x_2^8]^{600} x_2
[x_1^8,x_2^8]^{700} x_2^{-1} [x_1^8,x_2^8]^{800} x_2^{-1}.
\endmultline
$$
We define an auxiliary word $u(x_1, x_2)$ as follows:
$$
\multline u(x_1,x_2) = \,[x_1^{48},x_2^{40}]^{100} x_1^6
[x_1^{48},x_2^{40}]^{200}
x_1^6 [x_1^{48},x_2^{40}]^{300} x_1^6 [x_1^{48},x_2^{40}]^{400} x_1^6 \\
[x_1^{48},x_2^{40}]^{500} x_2^5 [x_1^{48},x_2^{40}]^{600} x_2^5
[x_1^{48},x_2^{40}]^{700} x_2^5 [x_1^{48},x_2^{40}]^{800} x_2^5.
\endmultline
\tag 1.1
$$
We then construct $v_n(x_1, \dots, x_n)$ as follows: If $n=2$ then
$$v_2(x_1, x_2) = w_2(x_1, x_2).
\tag 1.2
$$
If $n=3$ then $$v_3(x_1, x_2, x_3) = u \Bigl( u \bigl( v_2(x_1,
x_2), v_2(x_2, x_3) \bigr), \, u \bigl( v_2(x_2, x_3), v_2(x_3,
x_1) \bigr) \Bigr). \tag 1.3
$$ Inductively, for $n \ge 4$, define
$$
\multline v_n(x_1, \dots, x_n) = u \Bigl( u \bigl( v_{n-1}(x_1,
x_2, x_3, \dots,
x_{n-1}),\, v_{n-1}(x_{n-1}, x_2, x_3, \dots, x_{n-2}, x_n) \bigr), \\
u \bigl( v_{n-1}(x_{n-1}, x_2, x_3, \dots, x_{n-2}, x_n), \,
v_{n-1}(x_n, x_2, x_3 \dots, x_{n-2}, x_1) \bigr) \Bigr).
\endmultline
\tag 1.4
$$
For instance, $$v_4(x_1, x_2, x_3, x_4) = u \Bigl( u \bigl(
v_3(x_1, x_2, x_3), v_3(x_3, x_2, x_4) \bigr), \, u \bigl(
v_3(x_3, x_2, x_4), v_3(x_4, x_2, x_1) \bigr) \Bigr). \tag 1.5
$$

In Section 2, we establish several technical lemmas concerning
properties of Ivanov's word $v_2(x_1, x_2)$ and the auxiliary word
$u(x_1, x_2)$ which will be used throughout this paper. In
Sections 3--5, we prove that, for each $n \ge 3$, the word
$v_n(x_1, \dots, x_n)$ constructed above is indeed a $C$-test word
with the property in the statement of the Theorem (the case $n=2$ is
already proved in [2]). We first treat the case $n=3$ in Section
~3, and then proceed by simultaneous induction on $n$ with base
$n=4$ together with several necessary lemmas in Sections 4--5.

Once the Theorem is proved, our Corollary 1 that is an extended
version of [2, Corollary 1] to the case where $\psi$ is an
endomorphism of $F_m$ with non-cyclic image follows immediately,
as intended, by taking $u=v_m(x_1, \dots, x_m)$:

\proclaim {Corollary 1} There exists an element $u \in F_m$ such
that if $\phi$ is an endomorphism, $\psi$ is an endomorphism of
$F_m$ with non-cyclic image, and $\phi (u) = \psi (u)$, then
$\phi$ also has non-cyclic image, more precisely, $\phi = \tau _S
\circ \psi$, where $S \in F_m$ is such that $\<S, \psi(u)\>$ is
cyclic, and $\tau_S$ is the inner automorphism of $F_m$ defined
by means of $S$.
\endproclaim

In Corollary 2, we provide a positive solution to the Shpilrain's problem mentioned above:

\proclaim {Corollary 2} There exist two elements $u_1, u_2 \in
F_m$ such that any endomorphism $\psi$ of $F_m$ with non-cyclic
image is uniquely determined by $\psi (u_1)$, $\psi (u_2)$.
\endproclaim

The proof of Corollary 2 makes use of the words $v_m(x_1, \dots,
x_m)$ and $v_{m+1}(x_1, \dots, x_{m+1})$. Its detailed proof is
given in Section 6. The idea and the techniques used in [2] are
developed further in the present paper.

\heading 2. Preliminary Lemmas
\endheading

We begin this section by establishing some notation and
terminology. We let $X$, $Y$ (with or without subscript) be words
in $F_m$ throughout this paper. By $X = Y$ we denote the equality
in $F_m$ of words $X$ and $Y$, and by $X \equiv Y$ the
graphical (letter-by-letter) equality of words $X$ and $Y$. The
length of a word $X$ is denoted by $|X|$ (note
$|x_1{x_1}^{-1}|=2$). We say that a word $X$ is a {\it proper
power} if $X=Y^{\ell}$ for some $Y$ with $\ell >1$ and that a
word $A$ is {\it simple} if $A$ is non-empty, cyclically reduced,
and is not a proper power. If $A$ is simple, then an $A$-{\it
periodic} word is a subword of $A^k$ with some $k > 0$.

Now let us introduce several lemmas concerning properties of the
word $v_2(x_1, x_2)$ defined by $(1.2)$. For proofs of Lemmas
1--2, see [2].

\proclaim {Lemma 1 [2, Lemma 3]} If the subgroup $\<X_1, X_2\>$ of
$F_m$ is non-cyclic, then $v_2(X_1, X_2)$ is neither equal to the
empty word nor a proper power. If $\<X_1, X_2\>$ is cyclic, then
$v_2(X_1, X_2)=~1$.
\endproclaim

\proclaim {Lemma 2 [2, Lemma 4]} If the subgroup $\<X_1, X_2\>$ of
$F_m$ is non-cyclic and $v_2(X_1, X_2)=v_2(Y_1, Y_2)$, then there
exists a word $Z \in F_m$ such that
$$Y_1=ZX_1Z^{-1} \quad \text{and} \quad Y_2=ZX_2Z^{-1}.$$
\endproclaim

\proclaim {Lemma 3} If the subgroup $\<X_1, X_2\>$ of $F_m$ is
non-cyclic, then $v_2(X_1, X_2)^{-1} \neq v_2(Y_1, Y_2)$ for any
words $Y_1$, $Y_2$.
\endproclaim

\demo {Proof} By way of contradiction, suppose that $v_2(X_1,
X_2)^{-1}=v_2(Y_1, Y_2)$ for some words $Y_1$, $Y_2$. If $\<Y_1,
Y_2\>$ is cyclic, then it follows from Lemma 1 that $v_2(Y_1,
Y_2)=1$, so that $v_2(X_1, X_2)=1$, i.e., $\<X_1, X_2\>$ is
cyclic. This contradiction to the hypothesis of the lemma allows
us to assume that $\<Y_1, Y_2\>$ is non-cyclic. As in [2, Lemmas
1--4], let $W$ be a cyclically reduced word that is conjugate to
$v_2(X_1, X_2)$, and let $B$ be a simple word such that $[X_1^8,
X_2^8]$ is conjugate to $B$ (recall from [3, 4] that a commutator
$[A, B]$ of two words $A, \, B$ is not a proper power). Then
according to [2, Lemma ~2], $W$ has the form
$$W \equiv R_1T_1R_2T_2 \cdots R_8T_8,$$ where $R_i$ are
$B$-periodic words with $(i \cdot 100-14)|B| < |R_i| \leq i \cdot
100|B|$, $0 \leq |T_i| < 6|B|$, and $3488|B|<|W|<3648|B|$. The
same holds for $v_2(Y_1, Y_2)$ and attach the prime sign $'$ to
the notations for $v_2(Y_1, Y_2)$. Then the equality $v_2(X_1,
X_2)^{-1}=v_2(Y_1, Y_2)$ yields that $W^{-1}$ is a cyclic
permutation of $W'$, so that $|W^{-1}|=|W'|$. At this point,
apply the arguments in [2, Lemma 4] to $W^{-1}$, $W'$ to get $B^{-1} \equiv B'$ and that $R_i'$ $B'$-overlaps only with
$R_i^{-1}$ for each $i=1, 2, \dots, 8$. But this is impossible,
for if $R_i'$ $B'$-overlapped only with $R_i^{-1}$ then
$R_{i+1}'$ would have to $B'$-overlap only with $R_{i-1}^{-1}$
(indices modulo 8) by the order of indices in $W^{-1} \equiv
T_8^{-1}R_8^{-1} \cdots T_2^{-1}R_2^{-1}T_1^{-1}R_1^{-1}$ and $W'
\equiv R_1'T_1'R_2'T_2' \cdots R_8'T_8'$. \qed
\enddemo

We also establish several lemmas concerning properties of the
auxiliary word $u(x_1, x_2)$ defined by $(1.1)$.

\proclaim {Lemma 4} If the subgroup $\<X_1, X_2\>$ of $F_m$ is
non-cyclic, then $u(X_1, X_2)$ is neither equal to the empty word
nor a proper power. If $\<X_1, X_2\>$ is cyclic, then either
$u(X_1, X_2)$ is equal to the empty word provided
$X_1^6=X_2^{-5}$ or otherwise $u(X_1, X_2)$ is a proper power.
\endproclaim

\demo {Proof} The proof of the first part is similar to that of
[2, Lemma 3], and the second part is immediate from definition
$(1.1)$ of $u(x_1, x_2)$. \qed
\enddemo

\proclaim {Lemma 5} If the subgroup $\<X_1, X_2\>$ of $F_m$ is
non-cyclic and $u(X_1, X_2)=u(Y_1, Y_2)$, then there exists a
word $Z \in F_m$ such that
$$Y_1=ZX_1Z^{-1} \quad \text{and} \quad Y_2=ZX_2Z^{-1}.$$
\endproclaim

\demo {Proof} Applying the same arguments as in [2, Lemma 4] to
$u(X_1, X_2)$ and $u(Y_1, Y_2)$, we deduce that
$Y_1^6=ZX_1^6Z^{-1}$ and $Y_2^5=ZX_2^5Z^{-1}$ for some word $Z
\in F_m$. Since extraction of roots is unique in a free group, it
follows that $Y_1= ZX_1Z^{-1}$ and $Y_2= ZX_2Z^{-1}$, as
required. \qed
\enddemo

\proclaim {Lemma 6} If the subgroup $\<X_1, X_2\>$ of $F_m$ is
non-cyclic, then $u(X_1, X_2)^{-1} \neq u(Y_1, Y_2)$ for any
words $Y_1$, $Y_2$.
\endproclaim

\demo {Proof} Suppose to the contrary that $u(X_1,
X_2)^{-1}=u(Y_1, Y_2)$ for some words $Y_1$, $Y_2$. Then the
subgroup $\<Y_1, Y_2\>$ of $F_m$ is non-cyclic, for otherwise the
equality $u(X_1, X_2)^{-1}=u(Y_1, Y_2)$ would yield that $u(X_1,
X_2)=u(Y_1, Y_2)^{-1}=u(Y_1^{-1}, Y_2^{-1})$, contrary to Lemma
4. From here on, follow the proof of Lemma 3 to arrive at a
contradiction. \qed
\enddemo

\proclaim {Lemma 7} If the subgroup $\<X_1, X_2\>$ of $F_m$ is
non-cyclic, then, for each $i=1,\, 2$, the subgroup $\<X_i,
u(X_1, X_2)\>$ of $F_m$ is also non-cyclic.
\endproclaim

\demo {Proof} Suppose on the contrary that $\<X_i, u(X_1, X_2)\>$
is cyclic for some $i=1, \, 2$. Then by Lemma ~4 we have
$$X_i=u(X_1, X_2)^{\ell} \tag 2.1
$$
for some non-zero integer $\ell$. Let $U$ be a cyclically reduced
word that is conjugate to $u(X_1, X_2)$, and let $C$ be a simple
word such that $[X_1^{48}, X_2^{40}]$ is conjugate to $C$. By
$\overline{X}$, we denote a cyclically reduced word that is
conjugate to a word $X$ in $F_m$. Then by the same arguments as
in [2, Lemma ~1],
$$\text {max}(|\overline{X_1}|, |\overline{X_2}|) < 6|C|,$$
and also by the same arguments as in [2, Lemma 2], $U$ has the
form $$U \equiv Q_1S_1Q_2S_2 \cdots Q_8S_8,$$ where $Q_i$ are
$C$-periodic words with $(i \cdot 100-14)|C| < |Q_i| \leq i \cdot
100|C|$, $0 \leq |S_i| < 6|C|$ and $3488|C|<|U|<3648|C|$. This
yields that $$|\overline{X_i}| < 6|C| < 3488|C| < |U| =
|\overline{u(X_1, X_2)}| \leq |\overline{u(X_1, X_2)^{\ell}}|,$$
which contradicts equality $(2.1)$. \qed
\enddemo

In the following lemma which will be useful in Sections 4--6, the
word $v_n(x_1, \dots, x_n)$ with $n \ge 3$ is defined by
$(1.3)$--$(1.4)$.

\proclaim {Lemma 8} If both $v_2(X_1, X_2)$ and $v_n(Y_1, \dots,
Y_n)$ with $n \ge 3$ are neither equal to the empty word nor
proper powers, then the subgroup $\<v_2(X_1, X_2), v_n(Y_1,
\dots, Y_n)\>$ of $F_m$ is non-cyclic.
\endproclaim

\demo {Proof} Suppose on the contrary that $\<v_2(X_1, X_2),
v_n(Y_1, \dots, Y_n)\>$ is cyclic. It then follows from the
hypothesis of the lemma that $$\text {either} \ \ v_2(X_1,
X_2)=v_n(Y_1, \dots, Y_n) \ \ \text {or} \ \ v_2(X_1,
X_2)^{-1}=v_n(Y_1, \dots, Y_n).$$ This implies, by definition
$(1.3)$--$(1.4)$ of $v_n(x_1, \dots, x_n)$, the existence of
words $Z_1$, $Z_2$ in $F_m$ such that $$\text {either} \ \
v_2(X_1, X_2)=u(Z_1, Z_2) \ \ \text {or} \ \ v_2(X_1,
X_2)^{-1}=u(Z_1, Z_2).$$ But since a similar argument to that in
Lemma 3 shows that the latter equality cannot hold, the former
must hold. From here on, follow the arguments in [2, Lemma 4] to
obtain that $$B^{-\a}\tilde{X_1}B^{\a}=\tilde{Z_1^6} \quad \text
{and} \quad B^{-\a}\tilde{X_1}^{-1}B^{\a}=\tilde{Z_1^6},$$ where
$\tilde{X_i}$, $\tilde{Z_i}$ are conjugates of $X_i$, $Z_i$,
respectively, and $B$ is a simple word such that
$B=[\tilde{X_1}^8, \tilde{X_2}^8]$. This yields
$\tilde{X_1}^2=1$, i.e., $\tilde{X_1}=1$, contrary to the
hypothesis $v_2(X_1, X_2) \neq 1$. \qed
\enddemo

\heading 3. The case $n=3$
\endheading

In this section, we prove that $v_3(x_1, x_2, x_3)$ is a $C$-test
word with the additional property that \linebreak $v_3(X_1, X_2,
X_3)=1$ if and only if the subgroup $\<X_1, X_2, X_3\>$ of $F_m$
is cyclic. We begin with lemmas that play crucial roles in proving
this assertion.

\proclaim {Lemma 9} If $u \bigl( v_2(X_1, X_2), v_2(Y_1, Y_2)
\bigr) =1$, then $v_2(X_1, X_2)=1$ and $v_2(Y_1, Y_2)=1$.
\endproclaim

\demo{Proof} The hypothesis of the lemma implies by Lemma 4 that
$v_2(X_1, X_2)^6=v_2(Y_1, Y_2)^{-5}$; hence if one of $v_2(X_1,
X_2)$ and $v_2(Y_1, Y_2)$ is equal to the empty word, then so is
the other. So assume that $v_2(X_1, X_2) \neq 1$ and $v_2(Y_1,
Y_2)\neq 1$. Notice that the equality $v_2(X_1, X_2)^6=v_2(Y_1,
Y_2)^{-5}$ implies that $\<v_2(X_1, X_2), v_2(Y_1, Y_2)\>$ is
cyclic. Hence, in view of Lemmas 1 and 3, we have $v_2(X_1,
X_2)=v_2(Y_1, Y_2)$. This together with $v_2(X_1, X_2)^6=v_2(Y_1,
Y_2)^{-5}$ yields that $v_2(X_1, X_2)=v_2(Y_1, Y_2)=1$, contrary
to our assumption. \qed
\enddemo

\proclaim {Lemma 10} Suppose that the subgroup $\<X_1, X_2,
X_3\>$ of $F_m$ is non-cyclic. Then $v_3(X_1, X_2, X_3) \neq ~1$.
Furthermore, either $v_3(X_1, X_2, X_3)$ is not a proper power or
it has one of the following three forms:
$$\alignat 2
(A1)\ v_3(X_1, X_2, X_3)&=v_2(X_2, X_3)^{960} &&\quad \text {and} \quad X_1=1;\\
(A2)\ v_3(X_1, X_2, X_3)&=v_2(X_3, X_1)^{400} &&\quad \text {and} \quad X_2=1;\\
(A3)\ v_3(X_1, X_2, X_3)&=v_2(X_1, X_2)^{576} &&\quad \text {and}
\quad X_3=1.
\endalignat
$$
\endproclaim

\proclaim {Remark} In view of Lemma 1, $v_2(X_2, X_3)$, $v_2(X_3,
X_1)$, $v_2(X_1, X_2)$ in $(A1)$, $(A2)$, $(A3)$, respectively,
are neither equal to the empty word nor proper powers.
\endproclaim

\demo {Proof} Recall from $(1.3)$ that $$v_3(X_1, X_2, X_3)= u
\Bigl( u \bigl(v_2(X_1, X_2), v_2(X_2, X_3) \bigr), \, u \bigl(
v_2(X_2, X_3), v_2(X_3, X_1) \bigr) \Bigr).$$ In the case where
the subgroup $\<u \bigl( v_2(X_1, X_2), v_2(X_2, X_3) \bigr), \,
u \bigl( v_2(X_2, X_3), v_2(X_3, X_1) \bigr) \>$ of $F_m$ is
non-cyclic, the assertion that $v_3(X_1, X_2, X_3)$ is neither
equal to the empty word nor a proper power, as desired, follows
immediately from Lemma 4. So we only need to consider the case
where $$\text {the subgroup} \ \<u \bigl( v_2(X_1, X_2), v_2(X_2,
X_3) \bigr), \, u \bigl( v_2(X_2, X_3), v_2(X_3, X_1) \bigr) \> \
\text {is cyclic.} \tag 3.1
$$
Here, if $u \bigl( v_2(X_1, X_2), v_2(X_2, X_3) \bigr)=u \bigl(
v_2(X_2, X_3), v_2(X_3, X_1) \bigr) =1$, then Lemma 9 implies
that $v_2(X_1, X_2)=v_2(X_2, X_3)=v_2(X_3, X_1)=1$; hence, by
Lemma 1, $\<X_1, X_2\>$, $\<X_2, X_3\>$ and $\<X_3, X_1\>$ are
all cyclic. This yields that $\<X_1, X_2, X_3\>$ is cyclic,
contrary to the hypothesis of the lemma. Thus, at least one of
the words $u \bigl( v_2(X_1, X_2), v_2(X_2, X_3) \bigr)$ and $u
\bigl( v_2(X_2, X_3), v_2(X_3, X_1) \bigr)$ has to be not equal
to the empty word. We divide this situation into three cases.

\proclaim {Case I} $u \bigl( v_2(X_1, X_2), v_2(X_2, X_3) \bigr)
\neq 1$ and $u \bigl( v_2(X_2, X_3), v_2(X_3, X_1) \bigr) =1$.
\endproclaim

It follows from $u \bigl( v_2(X_2, X_3), v_2(X_3, X_1) \bigr)=1$
and Lemma 9 that $v_2(X_2, X_3)=v_2(X_3, X_1)=1$; so, by Lemma 1,
$\<X_2, X_3\>$ and $\<X_3, X_1\>$ are cyclic. Since $\<X_1, X_2,
X_3\>$ is non-cyclic, $X_3$ must be equal to the empty word;
hence we have
$$\split
v_3(X_1, X_2, X_3)
&= u \Bigl( u \bigl( v_2(X_1, X_2), 1 \bigr), \, 1 \Bigr)= u \bigl( v_2(X_1, X_2)^{24}, 1 \bigr) \\
&= v_2(X_1, X_2)^{24 \cdot 24} =v_2(X_1, X_2)^{576}.
\endsplit
$$
Therefore $v_3(X_1, X_2, X_3)$ has form $(A3)$ in this case.

\proclaim {Case II} $u \bigl( v_2(X_1, X_2), v_2(X_2, X_3) \bigr)
=1$ and $u \bigl( v_2(X_2, X_3), v_2(X_3, X_1) \bigr) \neq 1$.
\endproclaim

Since $u \bigl( v_2(X_1, X_2), v_2(X_2, X_3) \bigr) =1$, we have,
by Lemmas 1 and 9, that $\<X_1, X_2\>$ and $\<X_2, X_3\>$ are
cyclic, so that $X_2=1$; thus
$$\split
v_3(X_1, X_2, X_3)
&= u \Bigl( 1, \, u \bigl( 1, v_2(X_3, X_1) \bigr) \Bigr)= u \bigl(1, v_2(X_3, X_1)^{20} \bigr) \\
&= v_2(X_3, X_1)^{20 \cdot 20} =v_2(X_3, X_1)^{400}.
\endsplit
$$
Therefore $v_3(X_1, X_2, X_3)$ has form $(A2)$ in this case.

\proclaim {Case III} $u \bigl( v_2(X_1, X_2), v_2(X_2, X_3)
\bigr) \neq 1$ and $u \bigl( v_2(X_2, X_3), v_2(X_3, X_1) \bigr)
\neq 1$.
\endproclaim

In this case, we want to prove:

\proclaim {Claim} This case is reduced to the following two cases:
$$\split
&\text {\rm (i)} \ \<v_2(X_1, X_2), v_2(X_2, X_3), v_2(X_3, X_1)\> \ \text {is cyclic}; \\
&\text {\rm (ii)} \ \text {both} \ \<v_2(X_1, X_2), v_2(X_2, X_3)\> \ \text
{and} \ \<v_2(X_2, X_3), v_2(X_3, X_1)\> \ \text{are non-cyclic}.
\endsplit
$$
\endproclaim

\demo {Proof of the Claim} Assuming at least one of $\<v_2(X_1, X_2),
v_2(X_2, X_3)\>$ and $\<v_2(X_2, X_3), v_2(X_3, X_1)\>$ is
cyclic, we want to show that Case (i) occurs. Let us say that
$\<v_2(X_1, X_2), v_2(X_2, X_3)\>$ is cyclic (the case where
$\<v_2(X_2, X_3), v_2(X_3, X_1)\>$ is cyclic is analogous). If
$v_2(X_2, X_3)=1$, then $u \bigl( v_2(X_1, X_2), v_2(X_2, X_3)
\bigr)= v_2(X_1, X_2)^{24}$ and $u \bigl( v_2(X_2, X_3), v_2(X_3,
X_1) \bigr)=v_2(X_3, X_1)^{20}$. It then follows from $(3.1)$
that $\<v_2(X_1, X_2), v_2(X_3, X_1)\>$ is cyclic, which means
that Case (i) occurs. Now let $v_2(X_2, X_3) \neq 1$. Then
since $\<v_2(X_1, X_2), v_2(X_2, X_3)\>$ is cyclic, $(3.1)$
yields by Lemma 7 that $\<v_2(X_2, X_3), v_2(X_3, X_1)\>$ is also
cyclic; hence $\<v_2(X_1, X_2), v_2(X_2, X_3), v_2(X_3, X_1)\>$
is cyclic, that is, Case (i) occurs as well. \qed
\enddemo

Case (i) is again divided into subcases according to the number
of non-empty words among $v_2(X_1, X_2)$, $v_2(X_2, X_3)$ and
$v_2(X_3, X_1)$. Here, we note that if there exists only one
non-empty word, then it has to be $v_2(X_2, X_3)$, for otherwise
we would have a contradiction to the hypothesis of Case III.
Therefore, Case III is decomposed into the following six subcases.

\proclaim {Case III.1} $v_2(X_1, X_2)=v_2(X_3, X_1)=1$ and
$v_2(X_2, X_3) \neq 1$.
\endproclaim

In this case, it follows from Lemma 1 that $\<X_1, X_2\>$ and
$\<X_3, X_1\>$ are cyclic, so that $X_1=1$; hence we have
$$\split
v_3(X_1, X_2, X_3)
&= u \Bigl( u \bigl( 1, v_2(X_2, X_3) \bigr), \, u \bigl( v_2(X_2, X_3), 1 \bigr) \Bigr)= u \bigl( v_2(X_2, X_3)^{20}, v_2(X_2, X_3)^{24} \bigr) \\
&= v_2(X_2, X_3)^{20 \cdot 24 + 24 \cdot 20}=v_2(X_2, X_3)^{960}.
\endsplit
$$
Thus, $v_3(X_1, X_2, X_3)$ has form $(A1)$.

\proclaim {Case III.2} $v_2(X_1, X_2)=1$ and $\<1 \neq v_2(X_2,
X_3), 1 \neq v_2(X_3, X_1)\>$ is cyclic.
\endproclaim

Since $v_2(X_1, X_2)=1$, we have by Lemma 1 that $\<X_1, X_2\>$
is cyclic. Also since \linebreak $\<1 \neq v_2(X_2, X_3), 1 \neq
v_2(X_3, X_1)\>$ is cyclic, we have $1 \neq v_2(X_2, X_3) =
v_2(X_3, X_1)$ by Lemmas 1 and 3. Apply Lemma 2 to this equality:
there exists a word $S \in F_m$ such that
$$X_2=SX_3S^{-1} \quad \text {and} \quad X_3=SX_1S^{-1},$$
which yields that $S^{-1}X_2S=SX_1S^{-1}$, so that
$X_2=S^2X_1S^{-2}$. It then follows from $\<X_1, X_2\>$ being
cyclic that $\<S, X_1, X_2\>$ is cyclic. This together with the
equality $X_3=SX_1S^{-1}$ implies that $\<X_1, X_2, X_3\>$ is
cyclic, contrary to the hypothesis of the lemma. Therefore this
case cannot occur.

\proclaim {Case III.3} $v_2(X_2, X_3)=1$ and $\<1 \neq v_2(X_1,
X_2), 1 \neq v_2(X_3, X_1)\>$ is cyclic.
\endproclaim

Repeat a similar argument to that in Case III.2 to conclude that
this case cannot occur.

\proclaim {Case III.4} $v_2(X_3, X_1)=1$ and $\<1 \neq v_2(X_1,
X_2), 1 \neq v_2(X_2, X_3)\>$ is cyclic.
\endproclaim

Also repeat a similar argument to that in Case III.2 to conclude
that this case cannot occur.

\proclaim {Case III.5} $\<1 \neq v_2(X_1, X_2), 1 \neq v_2(X_2,
X_3), 1 \neq v_2(X_3, X_1) \>$ is cyclic.
\endproclaim

In this case, it follows from Lemmas 1 and 3 that
$$1 \neq v_2(X_1, X_2)=v_2(X_2, X_3)=v_2(X_3, X_1).
\tag 3.2
$$
Applying Lemma 2 to these equalities, we have the existence of
words $T_1$ and $T_2$ in $F_m$ such that
$$\aligned
X_1=T_1X_2T_1^{-1}, \quad X_2=T_1X_3T_1^{-1}; \\
X_2=T_2X_3T_2^{-1}, \quad X_3=T_2X_1T_2^{-1}.
\endaligned
\tag 3.3
$$
Combining $(3.2)$ and $(3.3)$ yields that $\<T_1^{-1}T_2, X_3\>$,
$\<T_1, v_2(X_2, X_3)\>$ and $\<T_2, v_2(X_3, X_1)\>$ are all
cyclic. At this point, we apply Ivanov's argument (see [2, pp.
403--404]) to obtain the following:

\proclaim {Claim (Ivanov)} $T_1=T_2$.
\endproclaim

\demo {Proof of the Claim} Suppose on the contrary that $T_1 \neq
T_2$. It then follows from $\<T_1^{-1}T_2, X_3\>$ being cyclic
that
$$X_3^{\ell_1}=(T_1^{-1}T_2)^{\ell_2}
\tag 3.4
$$
with nonzero integers $\ell_1$ and $\ell_2$. It also follows from
$\<T_1, v_2(X_2, X_3)\>$ and $\<T_2, v_2(X_3, X_1)\>$ being cyclic
that
$$T_1 = v_2(X_2, X_3)^{\ell_3} \quad \text {and} \quad T_2 = v_2(X_3, X_1)^{\ell_4},$$
with integers $l_3$ and $l_4$ at least one of which is non-zero,
so that
$$T_1^{-1}T_2 = v_2(X_2, X_3)^{-\ell_3}v_2(X_3, X_1)^{\ell_4}.$$ Hence, by $(3.4)$,
$$X_3^{\ell_1} = [v_2(X_2, X_3)^{-\ell_3}v_2(X_3, X_1)^{\ell_4}]^{\ell_2}.
\tag 3.5
$$

Note, by definition $(1.2)$ of $v_2(x_1, x_2)$, that the right
hand side of equality $(3.5)$ belongs to the subgroup $[F_m,
{\Cal N}_3]$, where ${\Cal N}_3$ is the normal closure in $F_m$
of the word $X_3$. So inside the relation module $\hat {\Cal N}_3
= {\Cal N}_3 / [{\Cal N}_3, {\Cal N}_3]$ of the one-relator group
$$G = \<x_1, \dots, x_m \| X_3\>,$$
equality $(3.5)$ can be expressed as
$$(\ell_1 - P) \cdot {\hat X}_3 = 0,$$
where $P$ is an element of the augmentation ideal of the group
ring ${\Bbb Z}(G)$ of $G$ over the integers and ${\hat X}_3$ is
the canonical generator of the relation module $\hat {\Cal N}_3$
of $G$. By Lyndon's result on the relation module $\Cal R$ of a
one-relator group $\<x_1, \dots, x_m \| R\>$ (see [4, 5]) which
says that if $Q \cdot \hat R = 0$ in $\Cal R$ then $Q$ is an
element of the augmentation ideal of ${\Bbb Z}(\<x_1, \dots, x_m
\| R\>)$, we must have $\ell_1=0$. This contradiction to the fact
$\ell_1 \neq 0$ completes the proof of the claim. \qed
\enddemo

If $T_1=T_2=1$, then equalities $(3.3)$ yield that $X_1=X_2=X_3$,
contrary to the hypothesis of the lemma. If $T_1=T_2 \neq 1$,
then we derive from equalities $(3.3)$ that
$$X_1=T_1^3X_1T_1^{-3},\ \ X_2= T_1^3X_2T_1^{-3} \ \ \text {and} \ \ X_3=T_1^3X_3T_1^{-3},$$
so that $\<T_1, X_1\>$, $\<T_1, X_2\>$ and $\<T_1, X_3\>$ are all
cyclic; therefore, $\<X_1, X_2, X_3\>$ is cyclic. A contradiction
implies that this case cannot occur.

\proclaim {Case III.6} Both $\<v_2(X_1, X_2), v_2(X_2, X_3)\>$
and $\<v_2(X_2, X_3), v_2(X_3\, X_1)\>$ are non-cyclic.
\endproclaim

In this case, in view of $(3.1)$ and Lemmas 4 and 6, we have that
$$u \bigl( v_2(X_1, X_2), v_2(X_2, X_3) \bigl) = u \bigl( v_2(X_2, X_3), v_2(X_3, X_1) \bigl),$$
so by Lemma 5 that there exists a word $T \in F_m$ such that
$$1 \neq v_2(X_1, X_2) = Tv_2(X_2, X_3)T^{-1} \ \ \text {and} \ \ 1 \neq v_2(X_2, X_3) = Tv_2(X_3, X_1)T^{-1}.
\tag 3.6
$$
Apply Lemma 2 to these equalities: there exist words $U_1$ and
$U_2$ in $F_m$ such that
$$\aligned
X_1=U_1X_2U_1^{-1}, \quad X_2=U_1X_3U_1^{-1}; \\
X_2=U_2X_3U_2^{-1}, \quad X_3=U_2X_1U_2^{-1}.
\endaligned
\tag 3.7
$$
Combining the equalities in $(3.6)$ and $(3.7)$, we deduce that
$\<U_1^{-1}U_2, X_3\>$, $\<T^{-1}U_1, v_2(X_2, X_3)\>$ and
$\<T^{-1}U_2, v_2(X_3, X_1)\>$ are cyclic. Here, apply Ivanov's
argument used in Case III.5 to get $U_1=U_2$. Then reasoning as
in Case III.5, we conclude that this case cannot occur.

The proof of Lemma 10 is complete. \qed
\enddemo

Now we are ready to prove the Theorem for the case $n=3$.

\demo {Proof of the Theorem ($n=3$)} The additional property that
$v_3(X_1, X_2, X_3)=1$ if and only if the subgroup $\<X_1, X_2,
X_3\>$ of $F_m$ is cyclic follows immediately from definition
$(1.3)$ of $v_3(x_1, x_2, x_3)$ and Lemma 10. So we only need to
prove that $v_3(x_1, x_2, x_3)$ is a $C$-test word, that is,
supposing $1 \neq v_3(X_1, X_2, X_3)=v_3(Y_1, Y_2, Y_3)$, we want
to prove the existence of a word $Z \in F_m$ such that
$$Y_i=ZX_iZ^{-1} \quad \text {for}\ \ \text {all} \ \ i=1,\, 2,\, 3.$$
We begin by dividing into two cases according to whether
$v_3(X_1, X_2, X_3)$ is a proper power or not.

\proclaim {Case I} $v_3(X_1, X_2, X_3)$ is a proper power.
\endproclaim

Applying Lemma 10 to $v_3(X_1, X_2, X_3)$ and $v_3(Y_1, Y_2,
Y_3)$, we have: $v_3(X_1, X_2, X_3)$ has one of three types
$(A1)$, $(A2)$ and $(A3)$; besides, by the equality $v_3(X_1, X_2,
X_3)=v_3(Y_1, Y_2, Y_3)$, $v_3(Y_1, Y_2, Y_3)$ has the same type
as $v_3(X_1, X_2, X_3)$, because the exponents in $(A1)$, $(A2)$
and $(A3)$ are all distinct. This gives us only three
possibilities $(A1) \& (A1)$, $(A2) \& (A2)$ and $(A3) \& (A3)$
for the types of $v_3(X_1, X_2, X_3) \, \& \, v_3(Y_1, Y_2, Y_3)$.

If $v_3(X_1, X_2, X_3)\, \& \, v_3(Y_1, Y_2, Y_3)$ is of type
$(A1) \& (A1)$ ($(A2) \& (A2)$ or $(A3) \& (A3)$ is similar), then
$$X_1=Y_1=1 \quad \text {and} \quad 1 \neq v_2(X_2, X_3)^{960}=v_2(Y_2, Y_3)^{960}.$$
Applying Lemma 2 to the equality $1 \neq v_2(X_2, X_3)=v_2(Y_2,
Y_3)$, we have that two $2$-tuples $(X_2, X_3)$ and $(Y_2, Y_3)$
are conjugate in $F_m$, which together with $X_1=Y_1=1$ yields
that two $3$-tuples $(X_1, X_2, X_3)$ and $(Y_1, Y_2, Y_3)$ are
conjugate in $F_m$, as desired.

\proclaim {Case II} $v_3(X_1, X_2, X_3)$ is not a proper power.
\endproclaim

In this case, it follows from Lemma 4 that
$$\text {the subgroup} \ \<u \bigl( v_2(X_1, X_2), v_2(X_2, X_3) \bigr), \, u \bigl( v_2(X_2, X_3), v_2(X_3, X_1) \bigr) \> \ \text {of} \ F_m \ \text {is non-cyclic}.
\tag 3.8
$$
This enables us to apply Lemma 5 to the equality $v_3(X_1, X_2,
X_3)=v_3(Y_1, Y_2, Y_3)$: for some word $S \in F_m$, we have
$$\aligned
1 \neq u \bigl( v_2(X_1, X_2), v_2(X_2, X_3) \bigr) = S u \bigl( v_2(Y_1, Y_2), v_2(Y_2, Y_3) \bigr) S^{-1}, \\
1 \neq u \bigl( v_2(X_2, X_3), v_2(X_3, X_1) \bigr) = S u \bigl(
v_2(Y_2, Y_3), v_2(Y_3, Y_1) \bigr) S^{-1}.
\endaligned
\tag 3.9
$$
Here, we consider two subcases.

\proclaim {Case II.1} One of $\<v_2(X_1, X_2), v_2(X_2, X_3)\>$
and $\<v_2(X_2, X_3), v_2(X_3, X_1)\>$ is non-cyclic.
\endproclaim

Let us say that $\<v_2(X_1, X_2), v_2(X_2, X_3)\>$ is non-cyclic
(the case where $\<v_2(X_2, X_3), v_2(X_3, X_1)\>$ is non-cyclic
is analogous). Then, by Lemma 5, the first equality of $(3.9)$
implies the existence of a word $W \in F_m$ such that
$$1 \neq v_2(X_1, X_2)=Wv_2(Y_1, Y_2)W^{-1} \quad \text {and} \quad 1 \neq v_2(X_2, X_3)=Wv_2(Y_2, Y_3)W^{-1}.
\tag 3.10
$$
Apply Lemma 2 to these equalities: there exist words $T_1$ and
$T_2$ in $F_m$ such that
$$\aligned
X_1=T_1Y_1T_1^{-1}, \quad X_2=T_1Y_2T_1^{-1}; \\
X_2=T_2Y_2T_2^{-1}, \quad X_3=T_2Y_3T_2^{-1}.
\endaligned
\tag 3.11
$$
Now applying Ivanov's argument introduced in Case III.5 of Lemma
10 to equalities $(3.10)$--$(3.11)$, we get $T_1=T_2$, by which
equalities $(3.11)$ yield the desired result.

\proclaim {Case II.2} Both $\<v_2(X_1, X_2), v_2(X_2, X_3)\>$ and
$\<v_2(X_2, X_3), v_2(X_3, X_1)\>$ are cyclic.
\endproclaim

In this case, if $v_2(X_2, X_3) \neq 1$, then the subgroup
$\<v_2(X_1, X_2), v_2(X_2, X_3), v_2(X_3, X_1)\>$ would be
cyclic, contrary to $(3.8)$. So $v_2(X_2, X_3)$ must be equal to
the empty word. On the other hand, equalities $(3.9)$ imply by
Lemma 4 that both $\<v_2(Y_1, Y_2), v_2(Y_2, Y_3)\>$ and
$\<v_2(Y_2, Y_3), v_2(Y_3, Y_1)\>$ are also cyclic. Then, for the
same reason as $v_2(X_2, X_3)$, $v_2(Y_2, Y_3)$ also has to be
equal to the empty word.

Thus, it follows from $(3.9)$ that
$$1 \neq v_2(X_1, X_2)^{24}=Sv_2(Y_1, Y_2)^{24}S^{-1} \quad \text {and} \quad 1 \neq v_2(X_3, X_1)^{20}=Sv_2(Y_3, Y_1)^{20}S^{-1},$$
that is,
$$1 \neq v_2(X_1, X_2)=Sv_2(Y_1, Y_2)S^{-1} \quad \text {and} \quad 1 \neq v_2(X_3, X_1)=Sv_2(Y_3, Y_1)S^{-1},$$
which is a similar situation to $(3.10)$. So from here on, we can
follow the proof of Case II.1 to obtain the desired result.

The proof of the Theorem for the case $n=3$ is complete. \qed
\enddemo

\heading 4. The base $n=4$ of simultaneous induction
\endheading

In this section, we prove the base step $n=4$ of simultaneous
induction which we use in Lemmas 11--13 and the Theorem.

\proclaim {Lemma 11 (n=4)} If both $v_2(X_1, X_2)$ and $v_3(Y_1,
Y_2, Y_3)$ are neither equal to the empty word nor proper powers,
then $\<v_2(X_1, X_2), v_3(Y_1, Y_2, Y_3)\>$ is non-cyclic.
\endproclaim

\demo {Proof} This is a special case of Lemma 8. \qed
\enddemo

\proclaim {Lemma 12 (n=4)} If $u \bigl( v_3(X_1, X_2, X_3),
v_3(Y_1, Y_2, Y_3) \bigr) =1$, $v_3(X_1, X_2, X_3)=v_3(Y_1, Y_2,
Y_3)=~1$.
\endproclaim

\demo {Proof} By Lemma 4, the hypothesis of the lemma implies that
$$v_3(X_1, X_2, X_3)^6=v_3(Y_1, Y_2, Y_3)^{-5},
\tag 4.1
$$
so that if one of $v_3(X_1, X_2, X_3)$ and $v_3(Y_1, Y_2, Y_3)$
is equal to the empty word, then so is the other. Hence assume
that $v_3(X_1, X_2, X_3) \neq 1$ and $v_3(Y_1, Y_2, Y_3) \neq 1$.

If one of $v_3(X_1, X_2, X_3)$ and $v_3(Y_1, Y_2, Y_3)$ is a
proper power, then so is the other by $(4.1)$, because $6$ and
$5$ are relatively prime. Hence by Lemma 10 $v_3(X_1, X_2, X_3)
\, \& \, v_3(Y_1, Y_2, Y_3)$ has nine possible types, and we can
easily check that $6$ times any of $960$, $400$ and $576$ never
equals $5$ times any of these, which means that equality $(4.1)$
cannot hold in any case, a contradiction.

If neither $v_3(X_1, X_2, X_3)$ nor $v_3(Y_1, Y_2, Y_3)$ is a
proper power, then by Lemma 6 we have $v_3(X_1, X_2,
X_3)=v_3(Y_1, Y_2, Y_3)$, because $(4.1)$ implies that
$\<v_3(X_1, X_2, X_3), v_3(Y_1, Y_2, Y_3)\>$ is cyclic. This
equality together with $(4.1)$ yields that $v_3(X_1, X_2, X_3)=
v_3(Y_1, Y_2, Y_3)=1$, contrary to our assumption. This completes
the proof. \qed
\enddemo

\proclaim {Lemma 13 (n=4)} Suppose that $\<X_1, X_2, X_3, X_4\>$
is non-cyclic. Then $v_4(X_1, X_2, X_3, X_4) \neq 1$.
Furthermore, either $v_4(X_1, X_2, X_3, X_4)$ is not a proper
power or it has one of the following four forms:
$$\alignat 2
(B1)\ v_4(X_1, X_2, X_3, X_4)&=v_2(X_4, X_3)^{960 \cdot 400} &&\quad \text {and} \quad X_1=X_2=1; \\
(B2)\ v_4(X_1, X_2, X_3, X_4)&=v_2(X_1, X_4)^{{400}^2} &&\quad \text {and} \quad X_2=X_3=1; \\
(B3)\ v_4(X_1, X_2, X_3, X_4)&=v_2(X_3, X_1)^{576 \cdot 400} &&\quad \text {and} \quad X_2=X_4=1; \\
(B4)\ v_4(X_1, X_2, X_3, X_4)&=v_3(X_1, X_2, X_3)^{1936} && \quad
\text {and} \quad X_1=X_3=X_4 \neq 1.
\endalignat
$$
\endproclaim

\proclaim {Remark} In view of Lemmas 1 and 10, $v_2(X_4, X_3)$,
$v_2(X_1, X_4)$, $v_2(X_3, X_1)$ and $v_3(X_1, X_2, X_3)$ in
$(B1)$, $(B2)$, $(B3)$ and $(B4)$, respectively, are neither
equal to the empty word nor proper powers.
\endproclaim

\demo {Proof} Recall from $(1.5)$ that
$$v_4(X_1, X_2, X_3, X_4) = u \Bigl( u \bigl( v_3(X_1, X_2, X_3), v_3(X_3, X_2, X_4) \bigr), \, u \bigl( v_3(X_3, X_2, X_4), v_3(X_4, X_2, X_1) \bigr) \Bigr).$$
If $\<u \bigl( v_3(X_1, X_2, X_3), v_3(X_3, X_2, X_4) \bigr), \,
u \bigl( v_3(X_3, X_2, X_4), v_3(X_4, X_2, X_1) \bigr)\>$ is
non-cyclic, then the assertion that $v_4(X_1, X_2, X_3, X_4)$ is
neither equal to the empty word nor a proper power, as desired,
follows directly from Lemma 4. So we only need to consider the
case where
$$\<u \bigl( v_3(X_1, X_2, X_3), v_3(X_3, X_2, X_4) \bigr), \, u \bigl( v_3(X_3, X_2, X_4), v_3(X_4, X_2, X_1) \bigr)\>\ \text {is cyclic}.
\tag 4.2
$$
Here, in order to avoid a contradiction to the hypothesis that
$\<X_1, X_2, X_3, X_4\>$ is non-cyclic, at least one of $u \bigl(
v_3(X_1, X_2, X_3), v_3(X_3, X_2, X_4) \bigr)$ and $u \bigl(
v_3(X_3, X_2, X_4), v_3(X_4, X_2, X_1) \bigr)$ has to be not
equal to the empty word. So we have three cases to consider.

\proclaim {Case I} $u \bigl( v_3(X_1, X_2, X_3), v_3(X_3, X_2,
X_4) \bigr) \neq 1$ and $u \bigl( v_3(X_3, X_2, X_4), v_3(X_4,
X_2, X_1) \bigr)=1.$
\endproclaim

In this case, we have, by Lemma 12 ($n=4$) and the Theorem ($n=3$),
that both $\<X_3, X_2, X_4\>$ and $\<X_4, X_2, X_1\>$ are cyclic.
Since $\<X_1, X_2, X_3, X_4\>$ is non-cyclic, $X_2$ and $X_4$
must be equal to the empty word; hence we have
$$\split
v_4(X_1, X_2, X_3, X_4)
&=u \Bigl( u \bigl( v_3(X_1, X_2, X_3), 1 \bigr), \, 1 \Bigr) \\
&=u \bigl( v_3(X_1, X_2, X_3)^{24}, 1 \bigr) = v_3(X_1, X_2, X_3)^{24 \cdot 24} \\
&=v_2(X_3, X_1)^{576 \cdot 400} \quad \text {by form}\ (A2).
\endsplit
$$
Therefore, $v_4(X_1, X_2, X_3, X_4)$ has form $(B3)$ in this case.

\proclaim {Case II} $u \bigl( v_3(X_1, X_2, X_3), v_3(X_3, X_2,
X_4) \bigr) = 1$ and $u \bigl( v_3(X_3, X_2, X_4), v_3(X_4, X_2,
X_1) \bigr) \neq 1.$
\endproclaim

In this case, by Lemma 12 ($n=4$) and the Theorem ($n=3$), we have
that both $\<X_1, X_2, X_3\>$ and $\<X_3, X_2, X_4\>$ are cyclic,
so that $X_2=X_3=1$; thus
$$\split
v_4(X_1, X_2, X_3, X_4)
&=u \Bigl( 1, \, u \bigl( 1, v_3(X_4, X_2, X_1) \bigr) \Bigr) \\
&=u \bigl( 1, v_3(X_4, X_2, X_1)^{20} \bigr) = v_3(X_4, X_2, X_1)^{20 \cdot 20} \\
&=v_2(X_1, X_4)^{400 \cdot 400} \quad \text {by form}\ (A2).
\endsplit
$$
Therefore, $v_4(X_1, X_2, X_3, X_4)$ has form $(B2)$ in this case.

\proclaim {Case III} $u \bigl( v_3(X_1, X_2, X_3), v_3(X_3, X_2,
X_4) \bigr) \neq 1$ and $u \bigl( v_3(X_3, X_2, X_4), v_3(X_4,
X_2, X_1) \bigr) \neq 1.$
\endproclaim

By reasoning as in Case III of Lemma 10, we break this case into the following six subcases.

\proclaim {Case III.1} $v_3(X_1, X_2, X_3)=v_3(X_4, X_2, X_1)=1$
and $v_3(X_3, X_2, X_4) \neq 1$.
\endproclaim

It follows from the Theorem ($n=3$) that $\<X_1, X_2, X_3\>$ and
$\<X_4, X_2, X_1\>$ are cyclic, so that $X_1=X_2=1$; hence we have
$$\split
v_4(X_1, X_2, X_3, X_4)
&=u \Bigl( u \bigl( 1, v_3(X_3, X_2, X_4) \bigr), \, u \bigl( v_3(X_3, X_2, X_4), 1 \bigr) \Bigr) \\
&=u \bigl( v_3(X_3, X_2, X_4)^{20}, v_3(X_3, X_2, X_4)^{24} \bigr)=v_3(X_3, X_2, X_4)^{20 \cdot 24 + 24 \cdot 20} \\
&=v_2(X_4, X_3)^{960 \cdot 400} \quad \text {by form}\ (A2).
\endsplit
$$
Thus, $v_4(X_1, X_2, X_3, X_4)$ has form $(B1)$.

\proclaim {Case III.2} $v_3(X_1, X_2, X_3)=1$ and $\<1 \neq
v_3(X_3, X_2, X_4), 1 \neq v_3(X_4, X_2, X_1)\>$ is cyclic.
\endproclaim

Since $v_3(X_1, X_2, X_3)=1$, we have, by the Theorem ($n=3$), that
$$\<X_1, X_2, X_3\> \ \ \text {is} \ \ \text {cyclic}.
\tag 4.3
$$
Also since $\<1 \neq v_3(X_3, X_2, X_4), 1 \neq v_3(X_4, X_2,
X_1)\>$ is cyclic, in view of Lemmas 10 and 11 ($n=4$), this case
is reduced to the following two cases:
$$\split
&\text {\rm (i)}\ \text {both} \ v_3(X_3, X_2, X_4)\ \text {and} \ v_3(X_4, X_2, X_1) \ \text {are proper powers}; \\
&\text {\rm (ii)}\ \text {neither} \ v_3(X_3, X_2, X_4)\ \text {nor} \
v_3(X_4, X_2, X_1) \ \text {is a proper power.}
\endsplit
$$

Case (i) is divided further into subcases according to types of
$v_3(X_3, X_2, X_4) \, \& \, v_3(X_4, X_2, X_1)$ by Lemma 10. Of
nine possible types of $v_3(X_3, X_2, X_4) \, \& \, v_3(X_4, X_2,
X_1)$, $(A3)\&(A1)$, $(A3)\&(A2)$ and $(A3)\&(A3)$ cannot occur,
for if $v_3(X_3, X_2, X_4)$ were of type $(A3)$, then $X_4=1$,
which together with $(4.3)$ yields a contradiction to the
hypothesis of the lemma. Also, $(A1)\&(A1)$ and $(A2)\&(A1)$
cannot occur, for if $v_3(X_4, X_2, X_1)$ were of type $(A1)$,
then $X_4=1$, again a contradiction. Moreover, $(A1)\&(A2)$
cannot occur, for this implies that $X_3=X_2=1$, so that
$v_3(X_3, X_2, X_4)=1$, a contradiction. Also, $(A2)\&(A3)$
cannot occur, for this implies that $X_2=X_1=1$, so that
$v_3(X_4, X_2, X_1)=1$, a contradiction as well.

For this reason, in Case (i), we only need to consider
$(A1)\&(A3)$ and $(A2)\&(A2)$ for the types of $v_3(X_3, X_2, X_4) \,
\& \, v_3(X_4, X_2, X_1)$. Therefore, Case III.2 is decomposed
into the following three subcases.

\proclaim {Case III.2.1} $v_3(X_3, X_2, X_4) \, \& \, v_3(X_4,
X_2, X_1)$ is of type $(A1) \& (A3)$.
\endproclaim

In this case, it follows from Lemma 10 that $X_3=X_1=1$,
$v_3(X_3, X_2, X_4)=v_2(X_2, X_4)^{960}$ and $v_3(X_4, X_2,
X_1)=v_2(X_4, X_2)^{576}$. Since $\<v_3(X_3, X_2, X_4), v_3(X_4,
X_2, X_1)\>$ is cyclic by the hypothesis of Case III.2, we have
that $\<v_2(X_2, X_4), v_2(X_4, X_2)\>$ is cyclic, so that $1
\neq v_2(X_2, X_4)=v_2(X_4, X_2)$ by Lemmas 1 and 3. Apply Lemma
2 to this equality: there is a word $S \in F_m$ such that
$$X_2=SX_4S^{-1} \quad \text {and} \quad X_4=SX_2S^{-1}.
\tag 4.4
$$
If $S=1$, then from $(4.4)$ we have $X_2=X_4$, which together
with $(4.3)$ yields a contradiction to the hypothesis of the
lemma. Now let $S \neq 1$. We derive from $(4.4)$ that
$$X_2=S^2X_2S^{-2} \quad \text {and} \quad X_4=S^2X_4S^{-2},$$
so that $\<S, X_2\>$ and $\<S, X_4\>$ are cyclic; thus $\<X_2,
X_4\>$ is cyclic. This together with $(4.3)$ yields a
contradiction as well (because $X_2 \neq 1$). Therefore, we
conclude that this case cannot occur.

\proclaim {Case III.2.2} $v_3(X_3, X_2, X_4)\, \& \, v_3(X_4,
X_2, X_1)$ is of type $(A2) \& (A2)$.
\endproclaim

It follows from Lemma 10 that $X_2=1$, $v_3(X_3, X_2,
X_4)=v_2(X_4, X_3)^{400}$ and $v_3(X_4, X_2, X_1)=v_2(X_1,
X_4)^{400}$. Since $\<v_3(X_3, X_2, X_4), v_3(X_4, X_2, X_1)\>$
is cyclic, $\<v_2(X_4, X_3), v_2(X_1, X_4)\>$ is cyclic, hence,
by Lemmas 1 and 3, $1 \neq v_2(X_4, X_3)=v_2(X_1, X_4)$. Then by
Lemma 2 there exists a word $U \in F_m$ such that
$$X_4=UX_1U^{-1} \quad \text {and} \quad X_3=UX_4U^{-1}.
\tag 4.5
$$
If $U=1$, then it follows from $(4.5)$ that $X_1=X_4=X_3$, which
together with $(4.3)$ yields a contradiction to the hypothesis of
the lemma. Now let $U \neq 1$. We have from $(4.5)$ that
$UX_1U^{-1}=U^{-1}X_3U$. This equality implies by $(4.3)$ that
$\<U, X_1, X_2, X_3\>$ is cyclic; thus, by the first equality of
$(4.5)$, we have that $\<X_1, X_2, X_3, X_4\>$ is cyclic. A
contradiction implies that this case cannot occur.

\proclaim {Case III.2.3} Neither $v_3(X_3, X_2, X_4)$ nor
$v_3(X_4, X_2, X_1)$ is a proper power.
\endproclaim

Since $\<v_3(X_3, X_2, X_4), v_3(X_4, X_2, X_1)\>$ is cyclic, we
have, by Lemma 6, that $$1 \neq v_3(X_3, X_2, X_4)=v_3(X_4, X_2,
X_1).$$ Apply the Theorem ($n=3$) to this equality: there is a word
$T \in F_m$ such that
$$X_3=TX_4T^{-1}, \ \ X_2=TX_2T^{-1} \ \ \text {and} \ \ X_4=TX_1T^{-1}.
\tag 4.6
$$
The second equality of $(4.6)$ implies that $\<T, X_2\>$ is
cyclic; hence, by $(4.3)$, $\<T, X_1, X_2, X_3\>$ is cyclic
(because $X_2 \neq 1$). Then by the third equality of $(4.6)$, we
have that $\<X_1, X_2, X_3, X_4\>$ is cyclic. A contradiction
implies that this case cannot occur.

\proclaim {Case III.3} $v_3(X_3, X_2, X_4)=1$ and $\<1 \neq
v_3(X_1, X_2, X_3), 1 \neq v_3(X_4, X_2, X_1)\>$ is cyclic.
\endproclaim

Repeat similar arguments to those in Case III.2 to conclude that
this case cannot occur.

\proclaim {Case III.4} $v_3(X_4, X_2, X_1)=1$ and $\<1 \neq
v_3(X_1, X_2, X_3), 1 \neq v_3(X_3, X_2, X_4)\>$ is cyclic.
\endproclaim

Also, repeat similar arguments to those in Case III.2 to conclude
that this case cannot occur.

\proclaim {Case III.5} $\<1 \neq v_3(X_1, X_2, X_3), 1 \neq
v_3(X_3, X_2, X_4), 1 \neq v_3(X_4, X_2, X_1)\>$ is cyclic.
\endproclaim

In this case, we want to prove:

\proclaim {Claim} $1 \neq v_3(X_1, X_2, X_3)=v_3(X_3, X_2,
X_4)=v_3(X_4, X_2, X_1).$
\endproclaim

\demo {Proof of the Claim} If none of these is a proper power, then
the assertion follows immediately from Lemma 6. So assume one of
these is a proper power. Then, in view of Lemmas 10 and 11
~($n=4$), the other two also have to be proper powers; thus two
of $X_1$, $X_3$ and $X_4$ must be equal to the empty word, unless
$X_2=1$. However, if two of $X_1$, $X_3$ and $X_4$ were equal to
the empty word, then we would have a contradiction to the
non-triviality of $v_3(X_1, X_2, X_3)$, $v_3(X_3, X_2, X_4)$ or
$v_3(X_4, X_2, X_1)$. Hence we must have $X_2=1$. Then
$$v_3(X_1, X_2, X_3)=v_2(X_3, X_1)^{400}, \ v_3(X_3, X_2, X_4)= v_2(X_4, X_3)^{400}, \ v_3(X_4, X_2, X_1)=v_2(X_1, X_4)^{400};
\tag 4.7
$$
hence the hypothesis of Case III.5 implies that $\<v_2(X_3, X_1),
v_2(X_4, X_3), v_2(X_1, X_4)\>$ is cyclic. It then follows from
Lemmas 1 and 3 that $v_2(X_3, X_1)=v_2(X_4, X_3)=v_2(X_1, X_4)$,
which together with $(4.7)$ proves the claim. \qed
\enddemo

Now apply the Theorem ($n=3$) to the equalities in the Claim: there exist
words $Z_1$ and $Z_2$ in $F_m$ such that
$$\aligned
X_1=Z_1X_3Z_1^{-1}, \quad X_2=Z_1X_2Z_1^{-1}, \quad X_3=Z_1X_4Z_1^{-1}; \\
X_3=Z_2X_4Z_2^{-1}, \quad X_2=Z_2X_2Z_2^{-1}, \quad
X_4=Z_2X_1Z_2^{-1}.
\endaligned
\tag 4.8
$$
We deduce from these equalities that $\<Z_1, X_2\>$, $\<Z_2,
X_2\>$, $\<Z_1Z_2^2, Z_1^2Z_2, X_1\>$, $\<Z_1Z_2^2, Z_1^2Z_2,
X_3\>$ and $\<Z_1Z_2^2, Z_1^2Z_2, X_4\>$ are all cyclic. Here, if
either $Z_1Z_2^2 \neq 1$ or $Z_1^2Z_2 \neq 1$, then we would have
that $\<X_1, X_2, X_3, X_4\>$ is cyclic, a contradiction. Hence
we must have that $Z_1Z_2^2=Z_1^2Z_2=1$, that is, $Z_1=Z_2=1$;
thus, by $(4.8)$,
$$X_1=X_3=X_4.$$

In addition, it follows from the Claim that
$$\split
v_4(X_1, X_2, X_3, X_4)
&=u \Bigl( u \bigl( v_3(X_1, X_2, X_3), v_3(X_1, X_2, X_3) \bigr), \, u \bigl( v_3(X_1, X_2, X_3), v_3(X_1, X_2, X_3) \bigr) \Bigr) \\
&=u \bigl( v_3(X_1, X_2, X_3)^{44}, v_3(X_1, X_2, X_3)^{44} \bigr) \\
&=v_3(X_1, X_2, X_3)^{44 \cdot 24 + 44 \cdot 20} \\
&=v_3(X_1, X_2, X_3)^{1936}.
\endsplit
$$
Therefore, in this case, $v_4(X_1, X_2, X_3, X_4)$ has form
$(B4)$.

\proclaim {Case III.6} Both $\<v_3(X_1, X_2, X_3), v_3(X_3, X_2,
X_4)\>$ and $\<v_3(X_3, X_2, X_4), v_3(X_4, X_2, X_1)\>$ are
non-cyclic.
\endproclaim

In view of $(4.2)$ and Lemmas 4 and 6, we have that
$$u \bigl( v_3(X_1, X_2, X_3), v_3(X_3, X_2, X_4) \bigr) = u \bigl( v_3(X_3, X_2, X_4), v_3(X_4, X_2, X_1) \bigr).$$
Then by Lemma 5 applied to this equality, there is a word $W \in
F_m$ such that
$$1 \neq v_3(X_1, X_2, X_3)=Wv_3(X_3, X_2, X_4)W^{-1}\ \ \text {and} \ \ 1 \neq v_3(X_3, X_2, X_4)=Wv_3(X_4, X_2, X_1)W^{-1}.$$
Applying the Theorem ($n=3$) to these equalities yields the existence
of words $V_1$ and $V_2$ in $F_m$ such that
$$\aligned
&X_1=V_1X_3V_1^{-1}, \quad X_2=V_1X_2V_1^{-1}, \quad X_3=V_1X_4V_1^{-1}; \\
&X_3=V_2X_4V_2^{-1}, \quad X_2=V_2X_2V_2^{-1}, \quad
X_4=V_2X_1V_2^{-1}.
\endaligned
$$
This is the same situation as $(4.8)$; hence, reasoning as in
Case III.5, we have $X_1=X_3=X_4$. But then $v_3(X_1, X_2,
X_3)=v_3(X_3, X_2, X_4)=v_3(X_4, X_2, X_1)$, which yields a
contradiction to the hypothesis of Case III.6. Therefore, we
conclude that this case cannot occur.

The proof of Lemma 13 ($n=4$) is now complete. \qed
\enddemo

\demo {Proof of the Theorem ($n=4$)} The additional property that
$v_4(X_1, X_2, X_3, X_4)=1$ if and only if the subgroup $\<X_1,
X_2, X_3, X_4\>$ of $F_m$ is cyclic is immediate from definition
$(1.5)$ of $v_4(x_1, x_2, x_3, x_4)$ and Lemma 13 ($n=4$). Now we
want to prove that $v_4(x_1, x_2, x_3, x_4)$ is a $C$-test word,
that is, supposing $1 \neq v_4(X_1, X_2, X_3, X_4)=v_4(Y_1, Y_2,
Y_3, Y_4)$, we want to prove the existence of a word $Z$ in $F_m$
such that
$$Y_i=ZX_iZ^{-1} \quad \text {for}\ \ \text {all} \ \ i=1,\, 2,\, 3, \, 4.$$
We consider two cases corresponding to whether $v_4(X_1, X_2,
X_3, X_4)$ is a proper power or not.

\proclaim {Case I} $v_4(X_1, X_2, X_3, X_4)$ is a proper power.
\endproclaim

Apply Lemma 13 ($n=4$) to $v_4(X_1, \dots, X_4)$ and $v_4(Y_1,
\dots, Y_4)$: $v_4(X_1, \dots, X_4)$ has one of the four types
$(B1)$--$(B4)$; besides, by the equality $v_4(X_1, \dots,
X_4)=v_4(Y_1, \dots, Y_4)$, $v_4(Y_1, \dots, Y_4)$ has the same
type as $v_4(X_1, \dots, X_4)$, since the exponents in
$(B1)$--$(B4)$ are all distinct. This gives us only four
possibilities $(B1)\&(B1), \cdots, (B4)\&(B4)$ for the types of
$v_4(X_1, \dots, X_4) \, \& \, v_4(Y_1,  \dots, Y_4)$.

If $v_4(X_1, \dots, X_4)\, \& \, v_4(Y_1, \dots, Y_4)$ is of type
$(B1)\&(B1)$ ($(B2)\&(B2)$ or $(B3)\&(B3)$ is analogous), then
$$X_1=X_2=Y_1=Y_2=1 \quad \text {and} \quad 1 \neq v_2(X_4, X_3)^{960 \cdot 400} = v_2(Y_4, Y_3)^{960 \cdot 400}.$$
Applying Lemma 2 to the equality $1 \neq v_2(X_4, X_3)=v_2(Y_4,
Y_3)$, we have that two 2-tuples $(X_4, X_3)$ and $(Y_4, Y_3)$
are conjugate in $F_m$, which together with $X_1=X_2=Y_1=Y_2=1$
yields the desired result.

If $v_4(X_1, \dots, X_4)\, \& \, v_4(Y_1, \dots, Y_4)$ is of type
$(B4)\&(B4)$, then
$$X_1=X_3=X_4, \ \ Y_1=Y_3=Y_4 \ \ \text {and} \ \ 1 \neq v_3(X_1, X_2, X_3)^{1936}=v_3(Y_1, Y_2, Y_3)^{1936}.$$
The equality $1 \neq v_3(X_1, X_2, X_3)=v_3(Y_1, Y_2, Y_3)$
yields, by the Theorem ($n=3$),  that two $3$-tuples $(X_1, X_2,
X_3)$ and $(Y_1, Y_2, Y_3)$ are conjugate in $F_m$. Then the
result follows from $X_1=X_3=X_4$ and $Y_1=Y_3=Y_4$.

\proclaim {Case II} $v_4(X_1, X_2, X_3, X_4)$ is not a proper
power.
\endproclaim

In view of Lemma 4, it follows that
$$\<u \bigl( v_3(X_1, X_2, X_3), v_3(X_3, X_2, X_4) \bigr), \, u \bigl( v_3(X_3, X_2, X_4), v_3(X_4, X_2, X_1) \> \ \text {is non-cyclic}.
\tag 4.9
$$
This enables us to apply Lemma 5 to the equality $v_4(X_1, \dots,
X_4)=v_4(Y_1, \dots, Y_4)$: there exists a word $S \in F_m$ such
that
$$\aligned
1 \neq u \bigl( v_3(X_1, X_2, X_3), v_3(X_3, X_2, X_4) \bigr) = S u \bigl( v_3(Y_1, Y_2, Y_3), v_3(Y_3, Y_2, Y_4) \bigr) S^{-1}, \\
1 \neq u \bigl( v_3(X_3, X_2, X_4), v_3(X_4, X_2, X_1) \bigr) = S
u \bigl( v_3(Y_3, Y_2, Y_4), v_3(Y_4, Y_2, Y_1) \bigr) S^{-1}.
\endaligned
\tag 4.10
$$
Here, we have four subcases to consider.

\proclaim {Case II.1} Both $\<v_3(X_1, X_2, X_3), v_3(X_3, X_2,
X_4)\>$ and $\<v_3(X_3, X_2, X_4), v_3(X_4, X_2, X_1)\>$ are
non-cyclic.
\endproclaim

The hypothesis of this case enables us to apply Lemma 5 to
equalities $(4.10)$: there exist words $T_1$ and $T_2$ in $F_m$
such that
$$\aligned
1 \neq v_3(X_1, X_2, X_3)=T_1v_3(Y_1, Y_2, Y_3)T_1^{-1}, \ \ 1 \neq v_3(X_3, X_2, X_4)=T_1v_3(Y_3, Y_2, Y_4)T_1^{-1}; \\
1 \neq v_3(X_3, X_2, X_4)=T_2v_3(Y_3, Y_2, Y_4)T_2^{-1}, \ \ 1
\neq v_3(X_4, X_2, X_1)=T_2v_3(Y_4, Y_2, Y_1)T_2^{-1}.
\endaligned
\tag 4.11
$$
Then by the Theorem ($n=3$) applied to $(4.11)$, there exist words
$U_1$, $U_2$ and $U_3$ such that
$$\aligned
X_1=U_1Y_1U_1^{-1}, \quad X_2=U_1Y_2U_1^{-1}, \quad X_3=U_1Y_3U_1^{-1}; \\
X_3=U_2Y_3U_2^{-1}, \quad X_2=U_2Y_2U_2^{-1}, \quad X_4=U_2Y_4U_2^{-1}; \\
X_4=U_3Y_4U_3^{-1}, \quad X_2=U_3Y_2U_3^{-1}, \quad
X_1=U_3Y_1U_3^{-1}.
\endaligned
\tag 4.12
$$
Here, if one of $U_1=U_2$, $U_2=U_3$ and $U_3=U_1$ is true, then
the required result follows directly from $(4.12)$. So assume
that $U_1$, $U_2$ and $U_3$ are pairwise distinct. Combining the
equalities in $(4.12)$, we deduce that $\<U_3^{-1}U_1, X_1\>$,
$\<U_1^{-1}U_2, X_2\>$, $\<U_2^{-1}U_3, X_2\>$, $\<U_3^{-1}U_1,
X_2\>$, $\<U_1^{-1}U_2, X_3\>$ and $\<U_2^{-1}U_3, X_4\>$ are all
cyclic, so that $\<X_1, X_2\>$, $\<X_3, X_2\>$ and $\<X_4, X_2\>$
are cyclic. Since $\<X_1, \dots, X_4\>$ is non-cyclic (this
follows from $v_4(X_1, \dots, X_4) \neq 1$), we must have
$X_2=1$; so, by $(4.12)$, $Y_2=1$. Then by Lemma 10, the
equalities on the first line of $(4.11)$ yield that
$$1 \neq v_2(X_3, X_1)^{400}=T_1v_2(Y_3, Y_1)^{400}T_1^{-1}, \ \ 1 \neq v_2(X_4, X_3)^{400}=T_1v_2(Y_4, Y_3)^{400}T_1^{-1},$$
namely,
$$1 \neq v_2(X_3, X_1)=T_1v_2(Y_3, Y_1)T_1^{-1}, \ \ 1 \neq v_2(X_4, X_3)=T_1v_2(Y_4, Y_3)T_1^{-1}.
\tag 4.13
$$
This is a similar situation to $(3.10)$, so from here on, we can
follow the proof of Case II.1 of the Theorem ~($n=3$) to obtain that
two 3-tuples $(X_1, X_3, X_4)$ and $(Y_1, Y_3, Y_4)$ are
conjugate in $F_m$. Since $X_2=Y_2=1$, the desired result follows.

\proclaim {Case II.2} $\<v_3(X_1, X_2, X_3), v_3(X_3, X_2,
X_4)\>$ is non-cyclic, and $\<v_3(X_3, X_2, X_4), v_3(X_4, X_2,
X_1)\>$ is cyclic.
\endproclaim

In this case, we can apply Lemma 5 to the first equality of
$(4.10)$: there exists a word $V \in F_m$ such that
$$1 \neq v_3(X_1, X_2, X_3)=Vv_3(Y_1, Y_2, Y_3)V^{-1}, \ \ 1 \neq v_3(X_3, X_2, X_4)=Vv_3(Y_3, Y_2, Y_4)V^{-1}.
\tag 4.14
$$
Then the Theorem ($n=3$) applied to $(4.14)$ yields the existence of
words $W_1$ and $W_2$ in $F_m$ such that
$$\aligned
X_1=W_1Y_1W_1^{-1}, \quad X_2=W_1Y_2W_1^{-1}, \quad X_3=W_1Y_3W_1^{-1}; \\
X_3=W_2Y_3W_2^{-1}, \quad X_2=W_2Y_2W_2^{-1}, \quad
X_4=W_2Y_4W_2^{-1}.
\endaligned
\tag 4.15
$$
If $W_1=W_2$, then the required result follows from $(4.15)$. Now
assume that $W_1 \neq W_2$. We deduce from $(4.15)$ that
$\<W_1^{-1}W_2, X_2\>$ and $\<W_1^{-1}W_2, X_3\>$ are cyclic, so
that
$$\<X_2, X_3\> \ \ \text {is} \ \ \text {cyclic.}
\tag 4.16
$$

On the other hand, since $\<v_3(X_3, X_2, X_4), v_3(X_4, X_2,
X_1)\>$ is cyclic by the hypothesis of Case ~II.2, in view of
Lemmas 10 and 11 ($n=4$), we have that $v_3(X_3, X_2, X_4)$ is a
proper power if and only if $v_3(X_4, X_2, X_1)$ is a proper
power. For this reason, this case is reduced to the following
three subcases.

\proclaim {Case II.2.1} $v_3(X_4, X_2, X_1)=1$ ($v_3(X_3, X_2,
X_4) \neq 1$ by the hypothesis of Case II.2).
\endproclaim

It follows from the Theorem ($n=3$) that $\<X_4, X_2, X_1\>$ is
cyclic. Since $\<X_2, X_3\>$ is cyclic by $(4.16)$, in order to
avoid a contradiction to the fact that $\<X_1, \dots, X_4\>$ is
non-cyclic, we must have
$$X_2=1;
\tag 4.17
$$
thus $Y_2=1$ by $(4.15)$. Then by Lemma 10, equalities $(4.14)$
yield that
$$1 \neq v_2(X_3, X_1)=Vv_2(Y_3, Y_1)V^{-1}, \ \ 1 \neq v_2(X_4, X_3)=Vv_2(Y_4, Y_3)V^{-1}.$$
This is the same situation as $(4.13)$; hence from here on,
repeating the proof of Case II.1, we obtain the desired result.

\proclaim {Case II.2.2} Both $v_3(X_3, X_2, X_4)$ and $v_3(X_4,
X_2, X_1)$ are proper powers.
\endproclaim

In view of Lemma 10, we have nine possibilities for the types of
$v_3(X_3, X_2, X_4) \, \& \, v_3(X_4, X_2, X_1)$. Of these nine
possible types, $(A1) \& (A1)$, $(A1) \&(A2)$, $(A1) \& (A3)$,
$(A2) \& (A1)$, $(A2) \& (A3)$, $(A3) \& (A2)$ and $(A3) \& (A3)$
cannot occur, for if one of these occurred, then two of $X_1$,
$X_2$, $X_3$ and $X_4$ should be equal to the empty word, which
yields a contradiction to the non-triviality of $v_3(X_1, X_2,
X_3)$, $v_3(X_3, X_2, X_4)$ or $v_3(X_4, X_2, X_1)$. So only $(A2)
\& (A2)$ and $(A3) \& (A1)$ can actually occur.

If $v_3(X_3, X_2, X_4) \, \& \, v_3(X_4, X_2, X_1)$ is of type
$(A2) \& (A2)$, then $X_2=1$, which is the same situation as
$(4.17)$. Hence from here on, following the proof of Case II.2.1,
we arrive at the desired result. If $v_3(X_3, X_2, X_4) \, \& \,
v_3(X_4, X_2, X_1)$ is of type $(A3) \& (A1)$, then $X_4=1$. The
result then follows from $(4.15)$.

\proclaim {Case II.2.3} Neither $v_3(X_3, X_2, X_4) \neq 1$ nor
$v_3(X_4, X_2, X_1) \neq 1$ is a proper power.
\endproclaim

Since $\<v_3(X_3, X_2, X_4), v_3(X_4, X_2, X_1)\>$ is cyclic, it
follows from Lemma 6 that $v_3(X_3, X_2, X_4)=v_3(X_4, X_2,
X_1)$, so from the second equality of $(4.10)$ that
$$1 \neq v_3(X_4, X_2, X_1)^{44}=S u \bigl( v_3(Y_3, Y_2, Y_4), v_3(Y_4, Y_2, Y_1) \bigr) S^{-1}.
\tag 4.18
$$
This equality implies by Lemma 4 that $\<v_3(Y_3, Y_2, Y_4),
v_3(Y_4, Y_2, Y_1)\>$ is also cyclic. We then observe that
equality $(4.18)$ can hold only when neither $v_3(Y_3, Y_2, Y_4)$
nor $v_3(Y_4, Y_2, Y_1)$ is a proper power and $v_3(Y_3, Y_2,
Y_4)=v_3(Y_4, Y_2, Y_1)$, by which $(4.18)$ yields that
$$1 \neq v_3(X_4, X_2, X_1)^{44}=Sv_3(Y_4, Y_2, Y_1)^{44}S^{-1},$$
namely,
$$1 \neq v_3(X_4, X_2, X_1)=Sv_3(Y_4, Y_2, Y_1)S^{-1}.$$
Now apply the Theorem ($n=3$) to this equality: there exists a word
$W_3 \in F_m$ such that
$$X_4=W_3Y_4W_3^{-1}, \quad X_2=W_3Y_2W_3^{-1}, \quad X_1=W_3Y_1W_3^{-1}.$$
Putting this together with $(4.15)$, we have the same situation
as $(4.12)$ except that we already assumed $W_1 \neq W_2$ in Case
II.2. Therefore, from here on, we can follow the proof of Case
II.1 to derive the result.

\proclaim {Case II.3} $\<v_3(X_1, X_2, X_3), v_3(X_3, X_2,
X_4)\>$ is cyclic, and $\<v_3(X_3, X_2, X_4), v_3(X_4, X_2,
X_1)\>$ is non-cyclic.
\endproclaim

It is sufficient to repeat similar arguments to those in Case
II.2 to arrive at the desired result.

\proclaim {Case II.4} Both $\<v_3(X_1, X_2, X_3), v_3(X_3, X_2,
X_4)\>$ and $\<v_3(X_3, X_2, X_4), v_3(X_4, X_2, X_1)\>$ are
cyclic.
\endproclaim

Arguing as in the proof of Case II.2 of the Theorem ($n=3$) replacing
$(3.8)$ and $(3.9)$ by $(4.9)$ and $(4.10)$, respectively, we
deduce that $v_3(X_3, X_2, X_4)=v_3(Y_3, Y_2, Y_4)=1$. So
$$\<X_3, X_2, X_4\> \ \ \text {is} \ \ \text {cyclic;}
\tag 4.19
$$
moreover, it follows from $(4.10)$ that
$$1 \neq v_3(X_1, X_2, X_3)^{24}=Sv_3(Y_1, Y_2, Y_3)^{24}S^{-1}, \ \ 1 \neq v_3(X_4, X_2, X_1)^{20}=Sv_3(Y_4, Y_2, Y_1)^{20}S^{-1},
$$
namely,
$$1 \neq v_3(X_1, X_2, X_3)=Sv_3(Y_1, Y_2, Y_3)S^{-1}, \ \ 1 \neq v_3(X_4, X_2, X_1)=Sv_3(Y_4, Y_2, Y_1)S^{-1}.
$$
Then by the Theorem ($n=3$) applied to these equalities, we have the
existence of words $Z_1$ and $Z_2$ in $F_m$ such that
$$\aligned
X_1=Z_1Y_1Z_1^{-1}, \quad X_2=Z_1Y_2Z_1^{-1}, \quad X_3=Z_1Y_3Z_1^{-1}; \\
X_4=Z_2Y_4Z_2^{-1}, \quad X_2=Z_2Y_2Z_2^{-1}, \quad
X_1=Z_2Y_1Z_2^{-1}.
\endaligned
\tag 4.20
$$
If $Z_1=Z_2$, then the result follows from $(4.20)$. Now assume
that $Z_1 \neq Z_2$. Then equalities $(4.20)$ yield that
$\<Z_1^{-1}Z_2, X_1\>$ and $\<Z_1^{-1}Z_2, X_2\>$ are cyclic, so
that $\<X_1, X_2\>$ is cyclic. Since $\<X_3, X_2, X_4\>$ is
cyclic by $(4.19)$, we must have $X_2=1$, which is the same
situation as $(4.17)$. Thus, from here on, we can follow the
proof of Case II.2.1 to get the required result.

The Theorem ($n=4$) is now completely proved. \qed
\enddemo

\heading 5. The inductive step
\endheading

In this section, we prove the inductive step of simultaneous
induction which we use in Lemmas 11--13 and the Theorem. Let $n \ge
5$ throughout this section.

\proclaim {Lemma 11} If both $v_{n-2}(X_1, \dots, X_{n-2})$ and
$v_{n-1}(Y_1, \dots, Y_{n-1})$ are neither equal to the empty
word nor proper powers, then $\<v_{n-2}(X_1, \dots, X_{n-2}),
v_{n-1}(Y_1, \dots, Y_{n-1})\>$ is non-cyclic.
\endproclaim

\demo {Proof} By way of contradiction, suppose that
$\<v_{n-2}(X_1, \dots, X_{n-2}), v_{n-1}(Y_1, \dots, Y_{n-1})\>$
is cyclic. Since both $v_{n-2}(X_1, \dots, X_{n-2})$ and
$v_{n-1}(Y_1, \dots, Y_{n-1})$ are non-proper powers, it follows
from Lemma 6 that
$$v_{n-2}(X_1, \dots, X_{n-2}) = v_{n-1}(Y_1, \dots, Y_{n-1}),
$$
so from $(1.3)$--$(1.4)$ and Lemmas 4 and 5 that there exists a
word $S \in F_m$ such that
$$
\aligned u \bigl( v_{n-3}(X_1, \dots, X_{n-3}), v_{n-3}(X_{n-3},
\dots, X_{n-2}) \bigr) &= S u \bigl( v_{n-2}(Y_1, \dots,
Y_{n-2}), v_{n-2}(Y_{n-2}, \dots, Y_{n-1}) \bigr) S^{-1},
\\
u \bigl( v_{n-3}(X_{n-3}, \dots, X_{n-2}), v_{n-3}(X_{n-2},
\dots, X_1) \bigr) &= S u \bigl( v_{n-2}(Y_{n-2}, \dots,
Y_{n-1}), v_{n-2}(Y_{n-1}, \dots, Y_1) \bigr) S^{-1}.
\endaligned
\tag 5.1
$$

We first assume that $\<v_{n-3}(X_1, \dots, X_{n-3}),
v_{n-3}(X_{n-3}, \dots, X_{n-2})\>$ is non-cyclic. This enables
us to apply Lemma 5 to the first equality of $(5.1)$: there
exists a word $T \in F_m$ such that
$$\aligned
v_{n-3}(X_1, \dots, X_{n-3})&=Tv_{n-2}(Y_1, \dots, Y_{n-2})T^{-1} \\
v_{n-3}(X_{n-3}, \dots, X_{n-2})&=Tv_{n-2}(Y_{n-2}, \dots,
Y_{n-1})T^{-1}.
\endaligned
\tag 5.2
$$
If both sides of the first equality of $(5.2)$ are non-proper
powers, then this equality yields a contradiction to the induction
hypothesis Lemma 11; if both sides of the first equality of
$(5.2)$ are proper powers, then we see from Lemma 10 and the
induction hypothesis of Lemma 13 that they cannot be the same
proper powers (because the exponents in both sides cannot be
identical), a contradiction as well.

We next assume that $\<v_{n-3}(X_1, \dots, X_{n-3}),
v_{n-3}(X_{n-3}, \dots, X_{n-2})\>$ is cyclic. Then the first
equality of $(5.1)$ yields by Lemma 4 that
$$\multline
\<v_{n-3}(X_1, \dots, X_{n-3}), v_{n-3}(X_{n-3}, \dots, X_{n-2}), \\
Sv_{n-2}(Y_1, \dots, Y_{n-2})S^{-1}, Sv_{n-2}(Y_{n-2}, \dots,
Y_{n-1})S^{-1}\> \ \ \text {is} \ \ \text{cyclic.}
\endmultline
\tag 5.3
$$
Here, in view of the induction hypothesis of Lemma 13, we see
that there are only two ways to avoid a contradiction to Lemma 8
and the induction hypothesis of Lemma 11: (i) any non-trivial
word of $Sv_{n-2}(Y_1, \dots, Y_{n-2})S^{-1}$ and
$Sv_{n-2}(Y_{n-2}, \dots, Y_{n-1})S^{-1}$ is of type $(B4)$;
(ii) any non-trivial word in $(5.3)$ is of one of types $(B1)$,
$(B2)$ and $(B3)$. (For $n=5$,  there is only one way: any
non-trivial word of $Sv_{n-2}(Y_1, \dots, Y_{n-2})S^{-1}$ and
$Sv_{n-2}(Y_{n-2}, \dots, Y_{n-1})S^{-1}$ is of one of types
$(A1)$, $(A2)$ and $(A3)$.) However, in either case, we can
observe that equalities $(5.2)$ cannot hold. This contradiction
completes the proof. \qed
\enddemo

\proclaim {Lemma 12} If $u \bigl( v_{n-1}(X_1, \dots, X_{n-1}),
v_{n-1}(Y_1, \dots, Y_{n-1}) \bigr)=1$, then $v_{n-1}(X_1, \dots,
X_{n-1})= 1$ and $v_{n-1}(Y_1, \dots, Y_{n-1})=1$.
\endproclaim

\demo {Proof} By Lemma 4, the hypothesis of the lemma yields that
$$v_{n-1}(X_1, \dots, X_{n-1})^6=v_{n-1}(Y_1, \dots, Y_{n-1})^{-5}.
\tag 5.4
$$
If one of $v_{n-1}(X_1, \dots, X_{n-1})$ and $v_{n-1}(Y_1, \dots,
Y_{n-1})$ is equal to the empty word, then by $(5.4)$ there is
nothing to prove. So assume that $v_{n-1}(X_1, \dots, X_{n-1})
\neq 1$ and $v_{n-1}(Y_1, \dots, Y_{n-1}) \neq ~1$. Since $6$ and
$5$ are relatively prime, equality $(5.4)$ implies that both
$v_{n-1}(X_1, \dots, X_{n-1})$ and $v_{n-1}(Y_1, \dots, Y_{n-1})$
are either proper powers or non-proper powers. If both are proper
powers, then a contradiction to equality $(5.4)$ follows from the
induction hypothesis of Lemma 13, for $6$ times any of $960 \cdot
400^{(n-4)}$, $400^{(n-3)}$, $576 \cdot 400^{(n-4)}$ and $1936$
cannot be equal to $5$ times any of these. If both are non-proper
powers, then, since $\<v_{n-1}(X_1, \dots, X_{n-1}), v_{n-1}(Y_1,
\dots, Y_{n-1})\>$ is cyclic, we have by Lemma 6 that
$v_{n-1}(X_1, \dots, X_{n-1})=v_{n-1}(Y_1, \dots, Y_{n-1})$. This
together with $(5.4)$ yields that $v_{n-1}(X_1, \dots,
X_{n-1})=v_{n-1}(Y_1, \dots, Y_{n-1})=1$. This contradiction to
our assumption completes the proof. \qed
\enddemo

\proclaim {Lemma 13} Suppose that $\<X_1, \dots, X_n\>$ is
non-cyclic. Then $v_n(X_1, \dots, X_n) \neq 1$. Furthermore,
either $v_n(X_1, \dots, X_n)$ is not a proper power or it has one
of the following four forms:
$$\aligned
(B1)\ v_n(X_1, \dots, X_n)&={\cases v_2(X_n, X_{n-1})^{960 \cdot
400^{(n-3)}}\quad \text {and} \quad
X_1=X_2=\cdots =X_{n-2}=1, &\text {for $n$ even} \\
v_2(X_{n-1}, X_n)^{960 \cdot 400^{(n-3)}}\quad \text {and} \quad
X_1=X_2=\cdots =X_{n-2}=1, &\text {for $n$ odd};
\endcases}
\\
(B2)\ v_n(X_1, \dots, X_n)&={\cases v_2(X_1,
X_n)^{400^{(n-2)}}\qquad \quad \ \text {and} \quad X_2=X_3=\cdots
=X_{n-1}=1, &\text {for $n$ even} \\
v_2(X_n, X_1)^{400^{(n-2)}}\qquad \quad \ \text {and} \quad
X_2=X_3=\cdots =X_{n-1}=1, &\text {for $n$ odd};
\endcases}
\\
(B3)\ v_n(X_1, \dots, X_n)&={\cases
v_2(X_{n-1}, X_1)^{576 \cdot 400^{(n-3)}}\quad \text {and} \quad X_2=\cdots =X_{n-2}=X_n=1, &\text  {for $n$ even} \\
v_2(X_1, X_{n-1})^{576 \cdot 400^{(n-3)}}\quad \text {and} \quad
X_2=\cdots =X_{n-2}=X_n=1, &\text {for $n$ odd};
\endcases}
\\
(B4)\ v_n(X_1, \dots, X_n)&=v_{n-1}(X_1, X_2, \dots,
X_{n-1})^{1936} \ \text {and} \quad X_1=X_{n-1}=X_n \neq 1.
\endaligned
$$
\endproclaim

\proclaim {Remark} In view of Lemma 1 and the induction
hypothesis of Lemma 13, $v_2(X_n, X_{n-1})$ and $v_2(X_{n-1},
X_n)$ in $(B1)$, $v_2(X_1, X_n)$ and $v_2(X_n, X_1)$ in $(B2)$,
$v_2(X_{n-1}, X_1)$ and $v_2(X_1, X_{n-1})$ in $(B3)$, and
$v_{n-1}(X_1, X_2, \dots, X_{n-1})$ in $(B4)$ are neither equal
to the empty word nor proper powers.
\endproclaim

\demo {Proof} Closely follow the proof of Lemma 13 ($n=4$)
replacing references to Lemmas 10, 11 ($n=4$), 12 ($n=4$) and
the Theorem ($n=3$) by references to the induction hypothesis of Lemma
13, Lemmas 8 and 11, Lemma 12 and the induction hypothesis of
the Theorem, respectively. Different situations from Lemma 13 ($n=4$)
can possibly occur only in Case III.2 and Case III.5, which we
reconsider below.

\proclaim {Case III.2} $\<1 \neq v_{n-1}(X_{n-1}, X_2, X_3,
\dots, X_{n-2}, X_n), \ 1 \neq v_{n-1}(X_n, X_2, X_3, \dots,
X_{n-2}, X_1)\>$ is cyclic, and $v_{n-1}(X_1, X_2, \dots,
X_{n-1})=1$.
\endproclaim

Since $v_{n-1}(X_1, X_2, \dots, X_{n-1})=1$, we have, by the
induction hypothesis of the Theorem, that
$$\<X_1, X_2, \dots, X_{n-1}\> \ \ \text {is} \ \ \text {cyclic}.
\tag 5.5
$$
Also since $\<1 \neq v_{n-1}(X_{n-1}, X_2, \dots, X_{n-2}, X_n),
1 \neq v_{n-1}(X_n, X_2, \dots, X_{n-2}, X_1)\>$ is cyclic, in
view of Lemmas 8 and 11 and the induction hypothesis of Lemma 13,
this case is reduced to the following two cases:
$$\split
&\text {\rm (i)}\ \text {both} \ v_{n-1}(X_{n-1}, X_2, \dots, X_{n-2}, X_n)\ \text {and} \ v_{n-1}(X_n, X_2, \dots, X_{n-2}, X_1) \ \text {are proper powers}; \\
&\text {\rm (ii)}\ \text {neither} \ v_{n-1}(X_{n-1}, X_2, \dots, X_{n-2},
X_n)\ \text {nor} \ v_{n-1}(X_n, X_2, \dots, X_{n-2}, X_1) \
\text {is a proper power.}
\endsplit
$$

Case (i) is divided further into subcases according to types of
$v_{n-1}(X_{n-1}, X_2, \dots, X_{n-2}, X_n)$ $\&$ $v_{n-1}(X_n,
X_2, \dots, X_{n-2}, X_1)$ by the induction hypothesis of Lemma
13. However, the former word cannot be of type $(B3)$ nor $(B4)$,
for this type together with $(5.5)$ yields a contradiction to the
hypothesis of Lemma 13. For the same reason, the latter cannot be
of type $(B1)$ nor $(B4)$. Also $(B1) \& (B2)$ and $(B2) \& (B3)$
cannot occur, for these types yield a contradiction to the
non-triviality of $v_{n-1}(X_{n-1}, X_2, \dots, X_{n-2}, X_n)$
and $v_{n-1}(X_n, X_2, \dots, X_{n-2}, X_1)$, respectively. Thus,
in Case (i), only $(B1) \& (B3)$ and $(B2) \& (B2)$ need to be
considered. This allows us to follow the proof of Lemma 13
($n=4$) from here on.

\proclaim {Case III.5} $\< \ 1 \neq v_{n-1}(X_1, X_2, \dots,
X_{n-1}), \ 1 \neq v_{n-1}(X_{n-1}, X_2, X_3, \dots, X_{n-2},
X_n),$ \linebreak $1 \neq v_{n-1}(X_n, X_2, X_3, \dots, X_{n-2},
X_1)\>$ is cyclic.
\endproclaim

In this case, in order to be able to keep following the proof of
Lemma 13 ($n=4$), it is sufficient prove the following:

\proclaim {Claim} $1 \neq v_{n-1}(X_1, X_2, \dots,
X_{n-1})=v_{n-1}(X_{n-1}, \dots, X_{n-2}, X_n)=v_{n-1}(X_n,
\dots, X_{n-2}, X_1)$.
\endproclaim

\demo {Proof of the Claim} If none of these is a proper power, then
the assertion follows immediately from Lemma 6. So assume that
one of these is a proper power. Then, in view of the induction
hypothesis of Lemma 13 and Lemmas 8 and 11, the other two also
have to be proper powers. Here, in order to avoid a contradiction
to the non-triviality of these words, all of these words must
have the same type either $(B2)$ or $(B4)$. The treatment of
$(B2)$ is the same as in the Theorem ($n=4$); if all of these words
are of type $(B4)$, then we have $X_1=X_{n-2}=X_{n-1}=X_n \neq
1$,  which proves the claim. \qed
\enddemo

The proof of Lemma 13 is complete. \qed
\enddemo

\demo {Proof of the Theorem} The additional property that $v_n(X_1,
\dots, X_n)=1$ if and only if the subgroup $\<X_1, \dots, X_n\>$
of $F_m$ is cyclic follows immediately from definition $(1.4)$ of
$v_n(x_1, \dots, x_n)$ and Lemma ~13. Now supposing $1 \neq
v_n(X_1, \dots, X_n)=v_n(Y_1, \dots, Y_n)$, we want to prove that
there exists a word $Z \in F_m$ such that
$$Y_i=ZX_iZ^{-1} \quad \text {for} \ \ \text {all} \ \ i=1, \, \dots, \, n,$$
that is, we want to show that $v_n(x_1, \dots, x_n)$ is a
$C$-test word. From here on, repeat the proof of the Theorem ($n=4$)
replacing references to Lemmas 10, 11 ($n=4$), 13 ($n=4$) and
the Theorem ~($n=3$) by references to Lemma 13, Lemmas 8 and 11, Lemma
13 and the induction hypothesis of the Theorem, respectively. A different
situation from the Theorem ($n=4$) can possibly occur only in Case
II.2.2, which we reconsider below.

\proclaim {Case II.2.2} Both $v_{n-1}(X_{n-1}, X_2, X_3, \dots,
X_{n-2}, X_n)$ and $v_{n-1}(X_n, X_2, X_3, \dots, X_{n-2}, X_1)$
are \linebreak proper powers.
\endproclaim

In view of Lemma 13, $v_{n-1}(X_{n-1}, X_2, \dots, X_{n-2}, X_n)$
$\&$ $v_{n-1}(X_n, X_2, \dots, X_{n-2}, X_1)$ has 16 possible
types, of which it suffices to consider the types involving
$(B4)$, namely, $(B1) \& (B4)$, $(B2) \& (B4)$, $(B3) \& (B4)$,
$(B4) \& (B4)$, $(B4) \& (B1)$, $(B4) \& (B2)$ and $(B4) \&
(B3)$, since the consideration for the remaining types is the
same as in the proof of the Theorem ($n=4$). However, none of these types except
for $(B4) \& (B4)$ can actually occur, for these types except for
$(B4) \& (B4)$ together with the hypothesis of Case II.2 that
$\<v_{n-1}(X_{n-1}, X_2, \dots, X_{n-2}, X_n), v_{n-1}(X_n, X_2,
\dots, X_{n-2}, X_1)\>$ is cyclic yield a contradiction to Lemma
~8.

If $(B4) \& (B4)$ occurs, then the equality corresponding to the
second one of $(4.10)$, namely,
$$\multline
1 \neq u \bigl( v_{n-1}(X_{n-1}, X_2, \dots, X_{n-2}, X_n), v_{n-1}(X_n, X_2, \dots, X_{n-2}, X_1) \bigr)\\
= S u \bigl( v_{n-1}(Y_{n-1}, Y_2, \dots, Y_{n-2}, Y_n),
v_{n-1}(Y_n, Y_2, \dots, Y_{n-2}, Y_1) \bigr) S^{-1}
\endmultline
$$
forces $v_{n-1}(Y_{n-1}, Y_2, \dots, Y_{n-2}, Y_n)\, \& \,
v_{n-1}(Y_n, Y_2, \dots, Y_{n-2}, Y_1)$ to have type $(B4) \&
(B4)$ as well. Thus, we have
$$X_1=X_{n-2}=X_{n-1}=X_n \neq 1 \quad \text {and} \quad Y_1=Y_{n-2}=Y_{n-1}=Y_n \neq 1,$$
which together with
$$X_{n-1}=W_2Y_{n-1}W_2^{-1}, \ \ X_2=W_2Y_2W_2^{-1}, \ \ \cdots, \ \ X_{n-2}=W_2Y_{n-2}W_2^{-1}, \ \ X_n=W_2Y_nW_2^{-1}$$
corresponding to the equalities on the second line of $(4.15)$
implies the desired result. This completes Case II.2.2. \qed
\enddemo

\heading 6. Proof of Corollary 2
\endheading

We are now in a position to prove Corollary 2.

\demo {Proof of Corollary 2} Let $\phi$ and $\psi$ be
endomorphisms of $F_m$ with non-cyclic images. Take
$$u_1=v_m(x_1, x_2, \dots, x_m) \quad \text {and} \quad u_2=v_{m+1}(x_{m-1}, x_{m-2}, \dots, x_1, x_m, x_1).$$
Supposing $\phi(u_i)=\psi(u_i)$, for $i=1, 2$, we want to prove
$\phi=\psi$. By Corollary 1, the equality $\phi(u_1)=\psi(u_1)$
implies that
$$\phi=\tau_S \circ \psi, \ \ \text {where\ \ $\<S, \psi(u_1)\>$\ \ is\ \ cyclic.}
\tag 6.1
$$
Here, it is sufficient to show $S=1$ to draw the desired result.
By $(6.1)$, the other equality $\phi(u_2)=\psi(u_2)$ yields that
$$\psi(u_2)=\phi(u_2)=S\psi(u_2)S^{-1},$$
so that $\<S, \psi(u_2)\>$ is cyclic. Then $S=1$ follows
obviously from the following:

\proclaim {Claim} $\<\psi(u_1), \psi(u_2)\>$ is non-cyclic.
\endproclaim

\demo {Proof of the Claim} Putting $X_i=\psi(x_i)$, for all
$i=1,2,\dots, m$, the claim is equivalent to
$$\<v_m(X_1, X_2, \dots, X_m), v_{m+1}(X_{m-1}, X_{m-2}, \dots, X_1, X_m, X_1)\> \ \ \text {is\ \ non-cyclic}.$$
Since $\psi$ has non-cyclic image, $\<X_1, X_2, \dots, X_m\>$ is
non-cyclic. This implies by the Theorem that
$$v_m(X_1, X_2, \dots, X_m) \neq 1 \ \ \text {and}\ \ v_{m+1}(X_{m-1}, X_{m-2}, \dots, X_1, X_m, X_1) \neq 1.$$
We treat two cases separately.

\proclaim {Case I} $v_{m+1}(X_{m-1}, X_{m-2}, \dots, X_1, X_m,
X_1)$ is not a proper power.
\endproclaim

In this case, $v_m(X_1, X_2, \dots, X_m)$ cannot be of form
$(B4)$, for if $v_m(X_1, X_2, \dots, X_m)$ were of form $(B4)$
then $X_1=X_{m-1}=X_m \neq 1$, which forces $v_{m+1}(X_{m-1},
X_{m-2}, \dots, X_1, X_m, X_1)$ to be a proper power of form
$(B4)$, a contradiction. Then the claim follows from Lemmas 8 and
11.

\proclaim {Case II} $v_{m+1}(X_{m-1}, X_{m-2}, \dots, X_1, X_m,
X_1)$ is a proper power.
\endproclaim

This case can occur only when $m \ge 3$, for if $v_3(X_1, X_2,
X_1)$ were a proper power, then we should have $X_2=1$, hence
$\<X_1, X_2\>$ is cyclic, contrary to our assumption that $\psi$
has non-cyclic image. It then follows from Lemma 13 that
$v_{m+1}(X_{m-1}, X_{m-2}, \dots, X_1, X_m, X_1)$ has one of four
types $(B1)$--$(B4)$. Of these types, $(B1)$ cannot occur, for
$(B1)$ implies $X_{m-1}=X_{m-2}=\cdots=X_1=1$, contrary to the
fact that $\<X_1, X_2, \dots, X_m\>$ is non-cyclic. For a similar
reason, $(B2)$ cannot occur, either. On the other hand, $(B4)$
cannot occur when $m=3$, for $(B4)$ together with $m=3$ implies
$X_2=X_3=X_1 \neq 1$, contrary to $\<X_1, X_2, X_3\>$ being
non-cyclic. When $m \ge 4$, $(B4)$ implies that $X_{m-1}=X_m=X_1
\neq 1$, so that $v_m(X_1, X_2, \dots, X_m)$ has type $(B4)$ as
well. Then the claim follows from Lemma 11.

It remains to consider type $(B3)$. If $v_{m+1}(X_{m-1}, X_{m-2},
\dots, X_1, X_m, X_1)$ has type $(B3)$, then we have that
$$\gather
v_{m+1}(X_{m-1}, X_{m-2}, \dots, X_1, X_m, X_1)= {\cases
v_2(X_{m-1}, X_m)^{576 \cdot 400^{(m-2)}}, &\text {for $m$ even} \\
v_2(X_m, X_{m-1})^{576 \cdot 400^{(m-2)}}, &\text {for $m$ odd},
\endcases}\\
\text {and} \quad X_{m-2}=X_{m-3}=\cdots=X_1=1. \tag 6.2
\endgather
$$
Then $(6.2)$ forces $v_m(X_1, X_2, \dots, X_m)$ to have type
$(B1)$; hence
$$v_m(X_1, X_2, \dots, X_m)=
{\cases
v_2(X_m, X_{m-1})^{960 \cdot 400^{(m-3)}}, &\text {for $m$ even} \\
v_2(X_{m-1}, X_m)^{960 \cdot 400^{(m-3)}}, &\text {for $m$ odd}.
\endcases}
$$
Now, by way of contradiction, suppose the contrary of the claim.
It then follows that
$$\<v_2(X_{m-1}, X_m), v_2(X_m, X_{m-1})\> \ \ \text {is\ \ cyclic,}$$
so by Lemmas 1 and 3 that
$$1 \neq v_2(X_{m-1}, X_m) = v_2(X_m, X_{m-1}).$$
Applying Lemma 2 to this equality yields the existence of $T \in
F_m$ such that
$$X_{m-1}=TX_mT^{-1} \quad \text {and} \quad X_m=TX_{m-1}T^{-1}.
\tag 6.3
$$
Combining these equalities, we see that both $\<X_{m-1}, T\>$ and
$\<X_m, T\>$ are cyclic. Here, if $T \neq 1$, then $\<X_{m-1},
X_m\>$ is cyclic; if $T=1$, then $X_{m-1}=X_m$ by $(6.3)$. This
together with $(6.2)$ yields a contradiction to the fact that
$\<X_1, X_2, \dots, X_m\>$ is non-cyclic, which completes the
proof of the claim. \qed
\enddemo

The proof of Corollary 2 is completed. \qed
\enddemo

\heading Acknowledgements
\endheading

The author would like to express her sincere thanks to Professor
S. V. Ivanov at the University of Illinois at Urbana-Champaign for his valuable comments and kind encouragement. She also thanks the referee for many helpful comments. 

\heading References
\endheading

\roster

\item"1." G. Baumslag, A. G. Myasnikov and V. Shpilrain, Open
problems in combinatorial group theory, {\it Contemp. Math.} {\bf
250} (1999), 1--27.

\item"2." S. V. Ivanov, On certain elements of free groups, {\it J. Algebra} {\bf 204} (1998), 394--405.

\item"3." A. G. Kurosh, ``The Theory of Groups, I, II,'' 2nd ed.,
Chelsea, New York, 1960.

\item"4." R. C. Lyndon and P. E. Schupp, ``Combinatorial Group Theory,'' Springer-Verlag, New York/Berlin, 1977.

\item"5." R. C. Lyndon, Cohomology theory of groups with a single defining relation, {\it Ann. of Math.} {\bf 52} (1962), 650--665.

\endroster

\enddocument